\documentclass[10pt,leqno,a4paper]{amsart}
\usepackage{amsmath}
\usepackage{amsthm}
\usepackage{amssymb}
\usepackage{graphicx}
\usepackage{epstopdf}
\usepackage{setspace}
\usepackage{hyperref}
\usepackage{pdfsync}
\usepackage{caption}
\usepackage{amsmath,amscd,amssymb}
\usepackage[english]{babel}
\usepackage[inner=3cm,outer=3cm,bottom=5cm]{geometry}
\usepackage{setspace}
\usepackage{mathrsfs}
\usepackage[all]{xy}
\usepackage{tikz-cd}

\newtheorem{theorem}{Theorem}[subsection]
\newtheorem{corollary}[theorem]{Corollary}
\newtheorem{lemma}[theorem]{Lemma}

\newtheorem{ass}[theorem]{Assumption}

\newtheorem{definition}[theorem]{Definition}
\newtheorem{conjecture}[theorem]{Conjecture}

\theoremstyle{remark}
\theoremstyle{definition}

\newtheorem{remark}[theorem]{Remark}
\newtheorem{example}[theorem]{Example}

\numberwithin{equation}{theorem}

\def\beq{\begin{equation}}
\def\eeq{\end{equation}}

\def\ben{\begin{enumerate}}
\def\een{\end{enumerate}}

\def\crash#1{}
\def\N{{\mathbb N}}
\def\Z{{\mathbb Z}}
\def\P{{\mathbb P}}

\def\R{{\mathbb R}}
\def\C{{\mathbb C}}
\def\A{{\mathbb A}}

\def\H{{\mathbb H}}

\def\l{\left}
\def\r{\right}
\def\[[{\l[\l[}
\def\]]{\r]\r]}

\def\spe{\text{spe}}
\def\sing{\text{Sing}}

\def\dg{{\dagger}}

\def\cf{\emph{cf.}\;}
\def\ie{\emph{i.e.}\;}
\def\lc{\emph{loc.cit.}\;}

\def\cA{{\mathcal A}}

\def\cF{{\mathcal F}}

\def\cO{{\mathcal O}}
\def\cH{{\mathcal H}}

\def\cS{{\mathcal S}}
\def\cT{{\mathcal T}}

\def\cU{{\mathcal U}}
\def\cX{{\mathcal X}}
\def\cY{{\mathcal Y}}

\def\sH{{\mathscr H}}

\def\fX{{\mathfrak X}}

\def\fc{{\mathfrak c}}

\def\wtilde{\widetilde}

\def\sm{{\rm sm}}

\def\Spf{{\rm Spf\,}}
\def\Spec{{\rm Spec\,}}

\def\iso{\xrightarrow{\ \sim\ }}

\def\char{{\text{char}}}

\def\limind{\mathop{\lim\limits_{\displaystyle\rightarrow}}}

\def\kc{{k^\circ}}
\def\kcc{{k^{\circ\circ}}}
\def\kt{\widetilde{k}}

\def\vphi{\varphi}
\def\eps{\epsilon}
\def\vt{{\vec{t}}}
\def\vv{{\vec{v}}}

\addtocounter{section}{-1}

\begin{document}

\setlength{\baselineskip}{0.55cm}	
\title{Riemann-Hurwitz formula for finite morphisms of $p$-adic curves}
\author{Velibor Bojkovi\'c}


\begin{abstract}
Given a finite morphism $\vphi:Y\to X$ of quasi-smooth Berkovich curves over a complete, non-archimedean, nontrivially valued algebraically closed field $k$ of characteristic 0, we prove a Riemann-Hurwitz formula relating their Euler-Poincar\'e characteristics (calculated using De Rham cohomology of their overconvergent structure sheaf). The main tools are $p$-adic Runge's theorem together with valuation polygons of analytic functions. Using the results obtained, we provide another point of view on Riemann-Hurwitz formula for finite morphisms of curves over algebraically closed fields of positive characteristic. 
\end{abstract}

\maketitle

\section{Introduction}

In algebraic, or complex analytic geometry, when studying finite morphisms between smooth proper curves, one of the most celebrated formulas that one encounters is the Riemann-Hurwitz formula (in further text RH formula), which relates the genera of the curves involved. That is, if $\vphi:Y\to X$ is a finite morphism between Riemann surfaces, RH formula states 
$$
\chi(Y)=\deg(\vphi)\cdot\chi(X)-\sum_{P\in Y}(e_p-1).
$$
Here, $\chi(Y)=2-2g(Y)$ and $\chi(X)=2-2g(X)$ are Euler-Poincar\'e characteristics and $g(Y)$ and $g(X)$ are genera of $Y$ and $X$, respectively, and $e_p$ is the ramification index (multiplicity) of the point $P$. A similar formula also holds for finite morphisms between smooth projective algebraic curves (over a field of characteristic 0). 

The purpose of this article is to provide an RH type of formula for a wide class of $p$-adic analytic curves in the sense of Berkovich. More precisely, we study finite morphisms $\vphi:Y\to X$ where $Y$ and $X$ are strict, connected quasi-smooth Berkovich $k$-analytic curves that admit a partition into finitely many affinoid domains with good canonical reduction and finitely many open annuli (\ie that admit finite triangulations, see Definition \ref{def:triangulation} and Assumption \ref{ass:triangulation}) and provide an RH formula which relates their Euler-Poincar\'e characteristics (see Theorem \ref{RH ft}). Here, $k$ is a complete, algebraically closed, non-archimedean and nontrivially valued field of characteristic 0. In particular, since projective smooth $k$-analytic curves admit finite triangulations and come by analytification of smooth projective $k$-algebraic curves, our formula may also be seen as a generalization of the classical RH formula for smoth projective $k$-algebraic curves.

The motivation as well as importance of studying finite morphisms between curves belonging to a wider class and not just proper ones and of having an RH formula for them comes from potential applications, especially  in the theory of $p$-adic differential equations (see for example \cite{Bal-Ked,Poi-Pul4}). In \lc authors have used special, simple cases of an RH formula for affinoid curves in order to prove the $p$-adic index theorem for curves of higher genus. It is very likely that such a formula will have a close relation with $p$-adic index theorem (either to be used to prove the index theorem in full generality \cf \cite{Bal-Ked,Poi-Pul4}, either to be a special case of the $p$-adic index theorem, as is shown in \cite{Boj16}). Let us mention that the RH formulas in \lc are just special cases of the formula that we present (Theorem \ref{RH ft}). 

We also emphasize the ubiquity of the formula in areas closely related to non-archimedean (Berkovich) geometry. For example, an RH formula appears in \cite[Corollary 7.3]{Hub01} for the case of adic curves (and is derived from the Lefschetz trace formula for $l$-adic sheaves), but also in \cite[Proposition 5.7]{Kat87} (proved using, among other tools, algebraic $K$-theory) for the case of two-dimensional Henselian normal local rings. 

Let us go back to our $p$-adic setting. As mentioned in the beginning, $Y$ and $X$ are $k$-analytic curves which admit a finite triangulation. We call such curves ft $k$-analytic curves (or just ft curves), and in particular quasi-smooth affinoid, smooth projective and wide open curves (but also many other) are ft curves. Finite morphisms between projective Berkovich curves and their RH formula (which in this case is the same as the algebraic one thanks to $p$-adic GAGA theorems) have been studied, most notably in \cite{JW,CTT14}. In \cite{JW}, the rich {\em inner} structure of Berkovich curves has been exploited to construct a pair of compatible deformation retractions of (projective) curves $Y$ and $X$ onto some of their respective skeleta, and furthermore using the classical RH formula for the corresponding algebraic curves, an RH type of formula is provided for such a compatible pair of skeleta. However, in \lc there is no discussion of a more general situation when the curves involved are not necessarily projective. 

The case of finite morphisms between quasi-smooth $k$-affinoid curves (but also wide open curves) and the corresponding RH formula appears in  \cite[Theorems 6.2.3 and 6.2.7]{CTT14}. Our Theorem \ref{RH ft} is a slight generalization of these formulas. We also point out that the results that coincide in \lc and in the present article, that is, RH formulas for quasi-smooth, compact and for wide open $k$-analytic curves, are done here independently, and using completely different methods. For example, the local case when $Y$ and $X$ are affinoid curves with good reduction in \lc is dealt with using algebraic geometry tools in the reduction of the morphism $\vphi$, while for us the main ingredient is the $p$-adic Runge's theorem combined with the classic RH formula. Again, passing from local to general case is different in \lc where the genus formula involving skeleta of the morphisms is exploited, while for us the main tool is the additivity of the Euler-Poincar\'e characteristics (Corollary \ref{EP add}) as well as the existence of the strictly $\vphi$-compatible triangulations (see Section 3 for definitions).

As we already mentioned, when $Y$ and $X$ are $k$-analytic projective curves, The RH formula is just the classic one. To see why it has to be modified in the case of a finite morphism of general ft  curves, let us discuss two very simple examples of finite morphisms between affinoid curves. At the same time we demonstrate the main ideas behind our proofs.

Suppose that the characteristic of the residual field is $\char(\kt)=p>0$ and suppose we are given two finite morphisms $f'_1:\P^1_k\to \P^1_k$, $x\mapsto 
y=x^p$, and $f'_2: \P^1_k\to \P^1_k$, $x\mapsto y=x^p-x$. 
Both maps induce finite morphisms $f_1, f_2: D(0,1^+)\to D(0,1^+)$ of the 
closed unit disc to itself.

The first map, when restricted to $D(0,1^+)$ is classically 
ramified\footnote{Classical ramification, \ie the ramification with 
support in rational points; classically ramified points are also called 
critical points, as in \cite{XF}.} only at $x=0$ with the 
ramification index $e_0=p$ (we distinguish the {\em finite} set of classically ramified points from the set of ramified points that consists of all the points where the morphism fails to be a local isomorphism. The locus of ramified points is much bigger than its subset of classically ramified points; in our example, the ramification locus 
contains all the points connecting $x=0$ with the Gauss point $\eta_{0,1}$
and even many more, see \cite{XF}), while the second one is classically unramified. 
Accepting for the moment that the 
Euler-Poincar\'e characteristic of the closed unit disc is equal to 1, the classical RH formula would give us $1=p\cdot 
1-\sum_{P\in D(0,1^+)(k)}(e_P-1)$ which is true only for the first map. 

To see what goes wrong with the second map, we could proceed as follows: deduce 
the RH formula for the map $f_2$ from the classical one by 
considering the full map $f_2'$ (fow which we know RH formula is true) and see how the relevant invariants change when we remove the open disc at infinity (\ie the complement of the closed unit disc). 

The classical ramification locus for the map $f_2'$ consists of the points 
$x_i=\zeta^i 
p^{-1/(p-1)}$ together with the point $x_{\infty}=\infty$, where $\zeta$ is any 
primitive $(p-1)$th root of 1 (and again, the full ramification locus is much 
richer). Each $x_i$ has the ramification index 
$e_{x_i}=2$, while $e_{x_{\infty}}=p$ so the classical Riemann-Hurwitz formula 
yields (recall $\chi(\P^1_k)=2$): $2=p\cdot 2-\sum_{i=1}^{p-1}(2-1)-(p-1)$  but 
for the moment let us write 
$\chi(\P^1_k)=\deg(f_2')\cdot\chi(\P^1)-\sum_{P\in D(\infty,1^-)(k)}(e_p-1)$, 
where $D(\infty,1^-)$ is the open unit disc with the center at $\infty$, i.e. 
$\P^1_k-  D(0,1^+)$. Using the fact that 
$\chi(D(0,1^+))=\chi(\P^1_k)-1$, we obtain 
$\chi(D(0,1^+))=\deg(f_2)\cdot\chi(D(0,1^+))-\big(\sum_{P\in 
D(\infty,1^-)(k)}(e_P-1)-\deg(f_2)+1\big)$, so we are led to think that the term 
$(\sum_{P\in 
D(\infty,1^-)(k)}(e_P-1)-\deg(f_2)+1)$ should count the defect. A natural question is: can 
we read off the defect out of the properties of the map $f_2$, without using 
its extension $f_2'$?

It turns out that the answer is positive under the condition that the map $f'_2$ extends to some extent over the boundary 
of the unit disc, and this leads us to stray in the area of finite \'etale 
morphisms of open annuli. Let us put here 
$A(0;r,1)$ to denote the open annulus $D(0,1^-)-  D(0,r^+)$ (note that 
every strict open annulus can be put in this form by a suitable isomorphism). 
Let $\vphi: 
A(0;r,1)\to A(0;r^n,1)$ be a finite \'etale morphism of degree $n$. After 
introducing suitable coordinates $T$ and $S$ at $0$, the derivative $dS/dT=d\vphi_{\#}(T)/d 
T$ is an invertible function on $A(0;r,1)$ and let us denote its 
order\footnote{Since $\frac{dS}{dT}$ is invertible, we can put it in the form 
$\frac{dS}{dT}=\eps T^{\sigma}(1+h(T))$, where for each $\rho \in (r,1)$, 
$|h(T)|_{\eta_{0,\rho}}<1$. Following the terminology of \cite[Lemma 1.6]{Lut93} we say that $\sigma$ is the order of 
$\frac{dS}{dT}$.} by 
$\sigma$. Suppose that our morphism $\vphi$ extends to a finite morphism $\psi: 
D(0,1^+)\to D(0,1^+)$, classically ramified at rational points $x_i$, 
$i=1,\dots, l$. Then, 
one can prove that $\sigma=\sum_{i=1}^l(e_{x_i}-1)$. 

Returning to our morphism $f_2'$, introducing coordinates at $\infty$ and 
considering the restriction of $f_2'$ to some annulus $A(0;r,1)$ 
(corresponding to the annulus $A(0;1,r')$ before introducing the 
coordinates at $\infty$, hence one can see the overconvergence of $f_2$ in
the fact that it extends also to some open annulus $A(0;1,r')$, for some 
$r'>1$), it follows that the aforementioned defect in fact is 
equal to the value $\nu:=\sigma-\deg(f_{2,\infty})+1$, where $f_{2,\infty}$ is 
the restriction of $f_2$ restricted to the open annulus $A(0;r,1)$ (with respect to new coordinates). 

Let us consider now a more general situation, a finite morphism $f:Y\to X$ of 
quasi-smooth, connected $k$-affinoid curves. By a result of Van Der Put 
\cite[Theorem 1.1]{VDP80} we know that $Y$ (resp $X$) can be 
embedded in a smooth,
projective curve $Y'$ (resp. $X'$), s.t. $Y'-  Y=\uplus_{i=1}^m D^Y_i$ 
(resp. $X'-  X=\uplus_{j=1}^nD^X_j$), where $D^Y_i$ (resp. $D^X_j$) are 
isomorphic to open unit discs. Suppose that $f$ extends to a finite morphism 
$f':Y'\to X'$. Let $f_i$ be the restriction of $f'$ to the disc $D^Y_i$, 
and let $\sigma_i$ be the order of its derivative when we restrict $f_i$ to a small 
open annulus living at the boundary of $D^Y_i$. Then, with similar arguments as 
before and with a bit more effort one can prove 
that 
\beq\label{intro RH}
\chi(Y)=\deg(f)\cdot\chi(X)-\sum_{i=1}^m\nu_i-\sum_{P\in Y(k)}(e_P-1)
\eeq
where $\nu_i=\sigma_i-\deg(f_i)+1$. This is a special case of our 
Riemann-Hurwitz formula \ref{RH} to which reader may now refer, or perhaps to \ref{RH ft} for the more general version. 

In the last example we made an essential assumption in order for the 
arguments to work, that the morphism $\vphi$ extends to a finite morphisms of {\em projective} 
curves $Y'$ and $X'$ which are made by adding {\em open discs} to the affinoids $Y$ and $X$, 
respectively. Although, thanks to Garuti's result \cite[Proposition 2.4]{Gar96}, it is always possible to do find a  {\em projectivization} of our morphism $\vphi$, that is to find $\vphi':Y'\to X'$ where so that $Y'$ and $X'$ are smooth and projective and contain $Y$ and $X$, respectively, so that $\vphi'_{|Y}=\vphi$, it is not always possible to find such $Y'$ by just adding open discs to $Y$. 

However, when $X$ is an affinoid subdomain in the projective line, our morphism $\vphi$ can be well {\em approximated} by a finite morphism $\vphi':Y'\to \P^1_k$, where $Y'-  Y$ is a disjoint union of open discs, and the approximation is due to the $p$-adic Runge's theorem. Our method is precisely this: to obtain the local RH formula, concerning affinoid curves with good reduction,  we simplify the situation and show that we can assume $X$ to be an affinoid domain in $\P^1_k$ and then deduce the RH formula for $\vphi$ from the classic one of the morphisms $\vphi':Y'\to \P^1_k$, by finding a good enough approximation $\vphi'$ of $\vphi$ (Definition \ref{def runge first}), and finally deducing the RH formula for $\vphi$ from the classic RH formula for $\vphi'$. Along the way of simplifying the general situation, we also prove the RH formula for the case where both $Y$ and $X$ are affinoid domains of $\P^1_k$ and here we use the full power of the theory of valuation polygons. Finally, the global case is done by using suitable, compatible partitions of $Y$ and $X$ into some simpler pieces, and deducing the global RH formula from the RH formulas for the restrictions of $\vphi$ over the elements of the partitions.

Up to now, we did not say what exactly we mean by EP characteristic of an ft curve. One has to face a choice which cohomology theory to use in order to get the right numbers, but having in mind future applications in the world of $p$-adic differential equations, we opted (out of few choices) for De Rham cohomology and we define EP characteristics as alternating sum of dimensions of De Rham cohomology groups calculated using the overconvergent structure sheaf of our curves. Although the results that we present in this article concerning De Rham cohomology and dimensions of De Rham cohomology groups of ft curves are certainly well known, we did not find a suitable reference (except in the case when the curves are affinoid or projective), hence the reason to collect such calculations and proofs in a separate section.\\

The article is structured as follows. 

The first section is mainly a remainder on the structure of quasi-smooth Berkovich curves which admit a finite triangulation. We characterize such curves (any such a curve is isomorphic to a complement of finitely many closed and finitely many open disjoint discs in a smooth projective curve), and we recall notions concerning them that we will need later on.

In the second section we study extensions/prolongations of finite morphisms. By a prolongation we simply mean a finite morphism defined on "bigger" curves which contain our starting curves $Y$ and $X$, and whose restriction on $Y$ is $\vphi$. We prove that finite morphisms of ft curves admit prolongations. Also, as was shown in the examples above, there are terms that are present in the RH formula but which are read off from the prolongation of the starting morphism (the terms $\nu_i$ in \eqref{intro RH}), and we need to show they do not depend on the chosen prolongation, and the proof of independence is postponed until the Section 5. 

The third section is concerned with finding triangulations of $Y$ and $X$ which behave nicely with respect to the morphism $\vphi$, what we call the strictly $\vphi$-compatible triangulations. We use them in order to reduce questions of RH formula for general ft curves to the case of curves which have good reduction.

The fourth section is dedicated to De Rham cohomology of ft curves, while the fifth one contains the discussion about RH formula, together with proofs. A separate subsection is dedicated to $p$-adic Runge's approximation theorem as well as to the discussion of how good approximations we actually need.

The final section is dedicated to the discussion of RH formula for smooth curves over the residue field $\wtilde{k}$. We do not distinguish whether the morphism is separable or purely inseparable.\\

\textbf{Acknowledgments}
 I would like to thank to my supervisor Francesco Baldassarri for proposing the problem of RH formula for finite morphisms of affinoid curves to me, and for his support during my work on this problem. I would also like to thank to my supervisor Denis Benois for his help during my stay in Bordeaux where the part of this work was done. I extend my thanks to Antoine Ducros, Marco Garuti, Elmar Gro{\ss}e-Kl{\"o}nne, J\'er\^ome Poineau and  Michael Temkin for answering my many questions and especially to J\'er\^ome and Antoine for their many suggestions on how to improve the present article. Finally I thank to the referee for many useful comments.

\section{Structure of quasi-smooth $k$-analytic curves}
\subsection{Base field} From now on, let $k$ be an algebraically closed field  which is complete with respect to a nonarchimedean and nontrivial valuation and of characteristic 0. The norm on $k$ will be denoted by $|\cdot|$, while $|\cdot|_{\infty}$ is reserved for the usual archimedean absolute value on $\R$. We denote by $\kc$ the set of integers of $k$, i.e. the set $\{a\in k, |a|\leq1\}$, and by $\kcc$ the maximal ideal of 
$\kc$, i.e. the set $\{a\in k, |a|<1\}$. The residual field $\kt$ is by definition 
$\kc/\kcc$ and it is an algebraically closed field.

\subsection{Basic pieces} Recall that the Berkovich projective 
line $\P^1_k$ is a one point compactification of the Berkovich affine line 
$\A^1_k$. The points of the affine line correspond to the multiplicative 
$k$-seminorms on the polynomial algebra $k[T]$. For $a\in k$ and $r\in 
\R_{\geq 0}$ we denote by $\eta_{a,r}$ the point in the affine line $\A^1_k$ 
corresponding to the multiplicative seminorm which is given by: for $f(t)\in 
k[T]$, $|f(\eta_{a,r})|=\max_{i\geq0}|\frac{f^{(i)}(a)}{i!}|r^i$. For the 
seminorm corresponding to the point $\eta_{a,\rho}$ we will write 
$|\cdot|_{a,\rho}$ or $|\cdot(\eta_{a,\rho})|$. If $a=0$ we also write 
$|\cdot|_\rho$ instead of $|\cdot|_{0,\rho}$. We identify 
points in $k$ with rational points in $\A^1_k$ via $a\in k \longleftrightarrow 
\eta_{a,0}\in\A^1_k$.

For $a\in k$ and $r\in\R_{\geq0}$ we denote by $D(a,r)$ (resp. $D(a,r^-)$) the 
 Berkovich closed (resp. open) disc centered in a point $a$ and of 
radius $r$. A point $\eta_{b,\rho}$ is in $D(a,r)$ (resp. $D(a,r^-)$) iff $b\in 
D(a,r)$ and $\rho\leq r$ (resp. $\rho<r$). Note that according to the previous, 
a rational point is a closed disc of radius 0.
Similarly, we denote by $A[a;r_1,r_2]$ (resp. $A(a;r_1,r_2)$) a closed (resp. 
open) Berkovich annulus with center in $a\in k$ and with radii $r_1$ and $r_2$, where $r_2\geq r_1\in \R_{>0}$ (resp. $r_2\geq 
r_1\in \R_{\geq 0}$). Note that the open annulus $A(a;0,r_2)$ is a punctured open disc \ie the 
open disc $D(a;r_2^-)$ punctured in the point $\eta_{a,0}$, which may differ from the classical notation where one considers $A(0;r_1,r_2)$ to be an open annulus only if $r_1>0$.

In general, let $X$ be a connected, quasi-smooth $k$-analytic curve. A subset 
$D\subset X$ is called a standard open disc (or just an open disc) if there 
is an isomorphism $T:D\iso D(0,r^-)$. A subset $A\in X$ is called a standard 
open annulus (or just an open annulus) if there is an isomorphism $T:A\iso 
A(0;r_1,r_2)$. A standard open disc $D$ is strict if there exists an isomorphism $D\iso D(0,1^-)$, and similarly, a standard open annulus is strict if there exists an isomorphism $A\iso A(0;,r_1,r_2)$, where both $r_1,r_2\in |k|$. 

Open discs and open annuli are special cases of {\em wide open curves} (or 
wide open spaces). We slightly generalize the definition given in  
\cite{Col03}.
\begin{definition}\label{wo} Let $X$ be a smooth projective curve (\ie analytification of a smooth projective $k$-algebraic curve). An open analytic subset $U$ of $X$ is called wide open (resp. strict wide open) in $X$ if 
$X-  U$ is 
a finite disjoint union (possibly empty) of closed discs (resp. strict closed discs) in $X$. In general, a $k$-analytic curve (\ie a 1-dimensional $k$-analytic space) is called a (strict) wide open curve if it is isomorphic to a (strict) wide open subset in some curve $X$ as before. 
\end{definition}

\subsection{Triangulations and semistable formal models}

Recall that if $X$ is a quasi-smooth curve and $U\subset X$ a subset, then the endpoints of $U$ in $X$ are points of the set $c(U)-U$, where we denote by $c(U)$ the closure of $U$ in $X$. Following \cite{Duc-book}, we denote by $X_{[2,3]}$ the set of type two and three points in a curve $X$ and we recall the definition of a triangulation which is slightly different from the one in \lc.
\begin{definition}\label{def:triangulation}
 Let $X$ be a quasi-smooth $k$-analytic curve. A closed and discrete subset $\cT\subset X_{[2,3]}$ is called a triangulation of $X$ if $X-\cT$ is a disjoint union of open discs and open annuli (\cf \cite[Section Triangulations]{Duc-book}). We say that $\cT$ is strict if every component in $X-\cT$ which is a relatively compact open annulus has exactly two endpoints. If for two triangulations $\cT$ and $\cT'$ of $X$ we have $\cT\subset \cT'$, we say that $\cT'$ is a refinement of $\cT$.
 \end{definition}
 \begin{remark}
  The difference with the definition in \lc is that we do not require the components of $X-\cT$ to be relatively compact in $X$. 
 \end{remark}

 We will use the following notation: for a quasi-smooth $k$-analytic curve $X$ and its triangulation $\cT$, $\cA_{\cT}(X)$ will denote the set of connected components in $X-\cT$ which are open annuli, while for $x\in \cT$, we denote by $W_{\cT,x}$, or just by $W_x$ if $\cT$ is understood from the context, the connected component of $X-\{\cT-\{x\}\}$ containing $x$. Note that if $x$ is not a point of the Shilov boundary of $X$, then $W_x$ is a wide open curve. If $\cT$ has at least two points, then for each point $x\in \cT$ we denote by $C_{\cT,x}$ the maximal affinoid domain in $X$ with Shilov boundary $x$ and such that $C_{\cT,x}\cap \cT=\{x\}$.
 
 \begin{ass}\label{ass:triangulation}
 It follows from Th\'eor\`eme 5.1.14 in \lc that if $X$ is strict then there exists a triangulation consisting just of type 2 points. Therefore we agree that whenever we speak about a triangulation of a strict curve, we assume that the triangulations consists only of type 2 points.
\end{ass}
When $X$ is a strict, connected, quasi-smooth compact $k$-analytic curve, we can assign a formal (strictly) semistable model of $X$ to any (strict) triangulation $\cT$ of $X$. The model, which we denote by $\fX_{\cT}$, is constructed as follows. We start by constructing formal schemes 
$\Spf\cO^{\circ}(W_\xi)$, for $\xi\in\cT$ and glue them 
along $\Spf\cO^{\circ}(W_{\xi_1}\cap W_{\xi_2})$ where $W_{\xi_1}\cap 
W_{\xi_2}\neq\emptyset$. The special fiber of 
$\fX_{\cT}$, denoted by $\fX_{\cT,s}$ is a reduced $\kt$-scheme, with only 
regular singular points. There is a well defined specialization map 
$\spe=\spe_{\fX_{\cT}}:X\to \fX_{\cT,s}$, some of whose properties are recalled in the following theorem  (\cite[Proposition 2.2 and 2.3]{Bo-Lu85} and \cite[Proposition 2.4.4]{Ber90}). 

\begin{theorem}\label{Bo-Lu-Be}{\em (Bosch-L\"utkebohmert-Berkovich)}
Let $X$ be a compact, connected, quasi-smooth $k$-analytic curve, $\cT$ a triangulation of 
$X$, $\fX_{\cT}$ the corresponding semistable formal model and 
$\spe:X\to\fX_{\cT,s}$ the specialization map. Let 
$x\in\fX_{\cT,s}$. Then
\begin{itemize}
 \item[(i)] The mapping $\spe$ induces a 1-1 correspondence between the 
irreducible components of $\fX_{\cT,s}$ (or generic points of irreducible 
components of $\fX_{\cT,s}$) and the points in $\cT$;

\item[(ii)] If $x$ is a smooth point in $\fX_{\cT,s}$ belonging to the 
irreducible component with generic point $\tilde{x}$, then $\spe^{-1}(x)$ is an 
open disc with the endpoint $\spe^{-1}(\tilde{x})$;
\item[(iii)] If $x$ is a regular singular point in $\fX_{\cT,s}$ belonging to 
the irreducible components with generic points $\tilde{x}_1$ and $\tilde{x}_2$ 
($\tilde{x}_1$ and $\tilde{x}_2$ may coincide), then $\spe^{-1}(x)$ is an 
open annulus with endpoints $\spe^{-1}(\tilde{x}_1)$ and 
$\spe^{-1}(\tilde{x}_2)$.
\end{itemize}
\end{theorem}

Using this theorem, given a (strictly) semistable model $\fX$ of $X$, we obtain the 
corresponding (strict) triangulation $\cT_{\fX}$ in the following 
way. Let $\{\fc_1,\dots,\fc_r\}$ be the irreducible components of $\fX_{s}$, 
let $\sm(\fc)$ denote the smooth part of the component $\fc$, and let 
$\sing(\fX_s)$ denote the singular locus of $\fX_s$ and let $\spe:X\to \fX_s$ be the specialization map with respect to the model $\fX$. If $\text{gen}(\fc_i)$ 
denotes the generic point of the component $\fc_i$, then 
$\cT_{\fX}=\{\spe^{-1}(\text{gen}(\fc_i)), i=1,\dots,r\}$. We 
furthermore have $A_{\cT_{\fX}}=\{\spe^{-1}(x),x\in\sing(\fX_s)\}$ and $W_{\cT_{\fX}}=\{\spe^{-1}(\fc_i),i=1,\dots,r\}$.

In the next section we will study triangulations in the context of finite morphisms.

\subsection{Ft curves}
Compact, connected quasi-smooth $k$-analytic curves have finite triangulations, but they are not the only curves with this property. As the curves that admit finite triangulations have the main role in this article, we name them in the next definition and characterize them in the theorem that follows.

\begin{definition}
 If $X$ is a strict, connected $k$-analytic curve, we say that it is an ft curve if $X$ admits a finite triangulation.
\end{definition}

\begin{theorem}\label{ft characterization}
If $X$ is a strict connected $k$-analytic curve which admits a finite triangulation, then $X$ is isomorphic to a complement of a disjoint union of finitely many strict open and closed discs in a smooth projective curve.
\end{theorem}
\proof
We may assume that $\cT$ is nonempty, since otherwise $X$ is isomorphic to a strict open disc or a strict open annulus (by definition), and in both cases the theorem is obviously true. First of all, because it admits a triangulation, every rational point is contained in a smooth analytic domain of $X$, hence $X$ is quasi-smooth. Second, for any such finite triangulation $\cT$, there are only finitely many open annuli in $\cA_{\cT}$. Indeed, since $X$ is a good $k$-analytic space \cite[Corollary 3.4]{DeJongEtale}, every point in $X$ has an affinoid neighborhood, and for $x\in \cT$ let $C_x$ be one such neighborhood. Lemma 3.1 in \lc implies that $C_x$ contains all but at most finitely many connected components in $X-\cT$ with an endpoint in $x$. In particular, there are at most finitely many open annuli that are connected components in $X-\cT$ with an endpoint in $x$ that are not contained in $C_x$, but also there are at most finitely many of them that are contained in $C_x$, the latter being easily verified. Since there are finitely many points in $\cT$, and each open annulus in $X-\cT$ must be attached to one of them, hence to intersect one of $C_x$, there are only finitely many elements in $\cA_{\cT}$.

Let us denote by $A_1,\dots,A_n$ the open annuli in $\cA_{\cT}$ which are not relatively compact in $X$. Note that for each $A_i$, the set $c(A_i)-A_i$ consists of only one point and it belongs to $\cT$. If we glue an open disc on the other side of the annulus $A_i$ we will obtain an open disc, which we glue to $X$ along $A_i$. We note that the glued open disc has an endpoint in $\cT$. In this way, we end up with a curve $X'$ which is quasi-smooth and such that $\cT$ is a triangulation of $X'$, but now there are no non-relatively compact annuli among the connected components in $X'-\cT$. We claim that $X'$ is compact. 

To prove the claim it is enough to show that $X'$ is a finite union of affinoid domains. Let $x\in \cT$ and let $C_x$ be an affinoid neighborhood of $x$. Once again we use that $C_x$ contains all but at most finitely many open discs in $X'-\cT$ which have an endpoint in $x$. These finitely many open discs can be covered with $C_x$ and finitely many sufficiently big closed discs contained in them. Furthermore, any open annulus in $X'-\cT$ can be covered by $C_x$, for $x\in \cT$ and a sufficiently big closed annulus contained in it. Hence the claim.

Now, $X'$ being compact implies further that it is either affinoid or projective curve \cite[Th\'eor\`eme 2]{F-M86}, and in the latter case we are done while in the former we may use \cite[Theorem 1.1]{VDP80} to conclude the proof of the theorem. 
\qed
\begin{remark} Consider an ft curve $X'$ and finitely many closed discs $D_i$, $i=1,\dots,n$ and finitely many open discs $E_j$, $j=1,\dots,m$ in $X'$ and let $X:=X'-\big(\cup_{i=1}^nD_i\cup_{j=1}^mE_j\big)$. Then, if $X'=\P^1_k$, $n=m=1$, $E_1=X'-D_1$, $X$ is an empty set and in particular not an ft curve. Contrary to this example, for every other choice for $X'$, $D_i$ and $E_j$, $X$ will be an ft curve.
 \end{remark}
\begin{definition}
If $X$ is an ft curve and $X'$ is a smooth projective $k$-analytic curve such that $X$ is isomorphic to a complement in $X'$ of a finite disjoint union of open and closed strict discs, then we say that $X'$ is a {\em simple projectivization} of $X$.
\end{definition}
\begin{ass}
Although some of the notions that we will recall later may hold in bigger generality, from now on we will assume that all the curves involved are {\em strict} ft curves, unless otherwise stated. 
\end{ass}

\subsection{Reduction} When $X$ is an affinoid curve the canonical reduction of $X$ is a $\tilde{k}$-algebraic curve $\wtilde{X}:=\Spec(A_{X}^\circ/A_{X}^{\circ\circ})$, where $A_X$ is the corresponding affinoid algebra of $X$ and $A^\circ_X:=\{f\in A_X| sp(f)\leq1\}$ and $A^{\circ\circ}_X:=\{f\in A_X| sp(f)<1\}$, where $sp$ is the spectral norm (see \cite[Section 2.4]{Ber90}). If $X$ is a quasi-smooth connected $k$-affinoid curve, then it has a good reduction if and only if (by definition) $\wtilde{X}$ is a smooth affine curve (and this happens if and only if $X$ has only one point in its Shilov boundary). 

\subsection{Skeleta} If $A$ is an annulus, recall that the skeleton of $A$, which we denote by $\Gamma^A$ is the complement of all open discs in $A$. For example, if $A=A(a;r_1,r_2)$, $\Gamma^A=\{\eta_{a,\rho},\rho\in (r_1,r_2)\}$. In general, for an ft curve $X$ and its triangulation $\cT$, we define the {\em skeleton} $\Gamma_{\cT}^X$ to be the union of the points in $\cT$ and all the skeleta of the open annuli in $\cA_{\cT}(X)$. If $X$ is clear from the context, we may just write $\Gamma_{\cT}$. In general, a skeleton of $X$ is any subset $\Gamma\subset X$ which is of the form $\Gamma=\Gamma^X_{\cT}$, for some triangulation $\cT$ of $X$. Similarly as before, for a point $x\in\Gamma_{\cT}-\cT$ we denote by $C_{\cT,x}$ the maximal affinoid domain in $X$ with the Shilov boundary equal to $x$ and having an empty intersection with $\cT$. In this case, $C_{\cT,x}$ is isomorphic to a closed annulus with equal inner and outer radii. \\

The following lemma is certainly well known, however as we will be using a similar reasoning in several places, we provide a complete proof.

\begin{lemma} \label{skeleta annuli}
 Let $\vphi:A_1\to A_2$ be a finite morphism of strict closed (resp. open) annuli. Then $\Gamma^{A_1}=\vphi^{-1}(\Gamma^{A_2})$.
\end{lemma}
\proof
Suppose that the annuli are closed, the other case being done in a similar fashion, and let $T:A[0;t,1]\iso A_1$ and $S:A[0;s,1]\iso A_2$ be some coordinates on $A_1$ and $A_2$, respectively. The morphism $\vphi$ can then be represented in coordinate form as an analytic function on $A_1$, $S=\vphi_{\#}(T)=\sum_{n\in \Z}a_n T^n=a_dT^d(1+h(T))$, where $|a_d|=1$ and $|(1+h(T))|_\rho<1$ for every $\rho\in [t,1]$ (see for example \cite[Lemma 1.6]{Lut93}). Here, $|d|_{\infty}$ is the degree of the morphism $\vphi$, and by changing the coordinate $S$ if necessary, we may assume that $d>0$. The skeleton of $A_1$ then is identified with the set of points $\eta_{0,\rho}$, where $\rho\in [t,1]$ and similarly, the skeleton of $A_2$ is identified with the set of points $\xi_{0,\rho'}$, $\rho'\in [s,1]$ (we changed the notation for the parameter on the target annulus from $\eta$ to $\xi$ in order to reduce the confusion). 
 
 Let us prove that $\vphi(\eta_{0,\rho})=\xi_{0,\rho^d}$, for $\rho\in [t,1]$. Let $\xi:=\vphi(\eta_{0,\rho})$. It is enough to prove that for every analytic function $f(S)$ on $A[0;s,1]$, $|f(\xi)|=|f(\xi_{0,\rho^d})|$. Moreover, since the polynomial functions are dense in the ring of analytic functions on $A[0;s,1]$ with respect to the spectral norm, it is enough to prove that for every $f(S)\in k[S]$, $|f(\xi)|=|f(\xi_{0,\rho^d})|$. Let us write $f(S)=b_0+\dots b_lS^l$. Then, 
 $$
 |f(\xi)|=|f(\vphi_{\#}(T))|_{0,\rho}=\max_{0\leq j\leq l}|b_j\Big(a_dT^d(1+h(T))\Big)^j|_{0,\rho}=\max_{0\leq j \leq l}|b_j||T^d|^j_{0,\rho}=|f(S)|_{0,\rho^d}.
 $$
 This proves that $\Gamma^{A_2}=\vphi(\Gamma^{A_1})$. Furthermore, since the valuation polygon of the function $\vphi_{\#}(T)$ has exactly one slope (and this slope is equal to $d$) we know that for every rational point $x\in A[0;t,1]$, the image of the disc $D(x,|x|^-)$ is contained in the annulus $A[0;|\vphi(x)|, |\vphi(x)|]$. Since $\vphi$ is continuous, this image is as well connected. We claim that $\vphi(D(x,|x|^-))=D(y,|y|^-)$, where $y=\vphi(x)$ (we may also note that $|y|=|x|^d$). Let $T_1$ and $S_1$ be coordinates on closed unit discs containing $A[0;t,1]$ and $A[0;s,1]$, respectively, so that $T_1(x)=0$ and $S_1(y)=0$ (hence, $D(x,|x|^-)$ is sent to $D(0,|x|^-)$ and $D(y,|y|^-)$ is sent to $D(0,|y|^-)$). The restriction of $\vphi$ on $D(x,|x|^-)$ expressed in new coordinates is a power series $S_1=\vphi_{\#1}(T_1)$ which is convergent on $D(0,|x|^-)$ and which sends $0$ to $0$. Once again, since the valuation polygon of the function $\vphi_{\#1}(T_1)$ is strictly increasing (by the choice of coordinates), convex and continuous function, if $z$ is any rational point in the image of $\vphi_{\#1}$, then by the theory of valuation polygons, every rational point in $D(0,|z|)$ is in the image as well. In particular, since rational points are dense in the respective sets, it follows that $\vphi(D(0,|x|^-))$ is an open disc. By continuity of $\vphi$, the set $\vphi\Big(D(0,|x|^-)\cup\{\eta_{0,|x|}\}\Big)=\vphi\Big(D(0,|x|^-)\Big)\cup\{\eta_{0,|y|}\}$ must be connected, and this proves the claim. 
 
 The fact that $\Gamma^{A_2}=\vphi(\Gamma^{A_1})$ and that connected components in $A_1-\Gamma^{A_1}$ are mapped to connected components in $A_2-\Gamma^{A_2}$ proves the lemma.
 \qed
 
 \begin{remark}  
 Using a similar reasoning with valuation polygons as in the previous lemma, one can prove that if $\vphi(T)$ is a holomorphic function on a closed annulus $A_1=A[0;r,1]$, then the image $\vphi(\Gamma^{A_1})=\Gamma^{A_2}$, where $A_2=A[0,r_1,r_2]$ and where $r_1=\min\{|\vphi(T)|_{\eta_{0,s}},s\in [r,1]\}$ and $r_2=\max\{|\vphi(T)|_{\eta_{0,r}},|\vphi(T)|_{\eta_{0,1}}\}$. A similar statement holds for an open annulus $A_1$ as well. 
\end{remark}

\subsection{Tangent space at a point} Let $X$ be an ft curve and $x\in X$ a point. 
The inductive limit over the open connected neighborhoods $U$ of $x$ in $X$ of the set of
connected components of $U- \{x\}$ is called the tangent space at $x$ and 
is denoted by $T_xX$. We will sometimes call the elements of the set $T_xX$ 
{\em the tangent points} and use bold letters $\vt$ and $\vv$ for them. If $x$ is of 
type 1 or 4, then the set $T_xX$ consists of one point \cite[Section 
4.2]{Bal-Ked}, this is because in these cases a fundamental system of 
neighborhoods can be chosen to consists of a decreasing sequence of discs $D_i$ in $X$ where 
$D_i- \{x\}$ has only one connected component. If $x$ is of type 3, 
then there exists an open neighborhood $A$ of $x$ in $X$ which is an open 
annulus with different endpoints. In this case, $T_xX$ consists of two elements, 
corresponding to the two connected components of $A- \{x\}$.

If $x$ is a type 2 point, the space $X$ can be described in another way. 
Namely, there exists a strict triangulation $\cT$ of $X$ that contains $x$. To show this, let us pick any strict triangulation $\cT'$ of $X$. If $x\in\Gamma_{\cT'}$, we can simply take $\cT:=\cT'\cup\{x\}$. Otherwise, $x$ belongs to an open disc $D$ which is a connected component in $X-\Gamma_{\cT'}$. If $D$ has an endpoint in $\cT'$, again we can take $\cT:=\cT'\cup\{x\}$, while if $D$ has an endpoint $\eta$ (necessarily of type 2) in $\Gamma_{\cT'}-\cT'$, we can take $\cT:=\cT'\cup\{x,\eta\}$. Next, from Theorem \ref{Bo-Lu-Be} it follows that $x$ corresponds to an irreducible component 
$\fc_x$ of $\fX_{\cT,s}$ and we note that if $U$ is a (small enough) connected open neighborhood of $x$, the connected components of $U- \{x\}$ are in 1-1 correspondence with the preimages of the rational points of $\fc_x(\kt)$ under
the specialization map from Theorem \ref{Bo-Lu-Be}. Hence, for type 2 point $x$ we may identify $T_xX$ with the rational points of $\fc_x$, or in other words, we may 
identify the tangent points in $T_xX$ with the residual classes which have 
an endpoint in $x$. 

It is useful to note that if $x$ is not in the interior of $X$, then $\fc_x$ is not projective. However, we may always embed $X$ in a smooth projective curve $X'$ and in this case the corresponding curve $\fc_x$ will be projective. The tangent space $T_xX'$ doesn't depend on the chosen curve $X'$, and the set of points $T_xX'-T_xX$ corresponds to the points at infinity of the smooth compactification of $\sm(\fc_x)$, see the remark below. We will denote the set $T_xX'$ by $T_xX^\dg$.  
  
Let $X'$ be a simple projectivization of $X$, so that $X'-X=\uplus_{1=1}^n B_i\uplus_{j=1}^m D_j$, where $B_i$ are open and $D_j$ are closed strict discs. For $i=1,\dots,n$, let $\xi_i\in X$ be the (necessarily) type two 
endpoint of the disc $B_i$, and let $\vt_i\in T_{\xi_i}X'$ be the tangent vector corresponding to $B_i$. For each $j=1,\dots,m$ let $\eta_j$ be the endpoint of $D_j$ and let $\vv_j\in T_{\eta_j}X'$ be the tangent vector "pointing" towards the curve $X$. Note that to each $\vv_j$ and for any triangulation $\cT$ of $X$, there exists a unique open annulus in $\cA_{\cT}(X)$ which has an endpoint in $X'-X$ and which corresponds to $\vv_j$.
  
\begin{definition}
 We define $T X$ to be the set of tangent vectors $\{\vt_1,\dots,\vt_n\}$ and $T_{in}X$ to be the set $\{\vv_1,\dots,\vv_m\}$, constructed in the previous paragraph. If $X$ compact we agree $T_{in}X$ to denote the empty set and if $X$ is  projective we also agree for $T X$ to denote the empty set. 
\end{definition}
\begin{remark}
 If $X$ is affinoid and $X'-  X=\uplus_{i=1}^n B_i $ as above, then the discs $B_i$ correspond to the points at infinity of the minimal compactification of the canonical reduction of $X$. 
\end{remark}

\section{Prolongations of finite morphisms}

Let $\vphi:Y\to X$ be a finite morphism of ft curves.
\begin{definition}
 A prolongation of the morphism $\vphi$ is a finite morphism of ft curves $\vphi' : Y' \to X'$,  together with open embeddings $i: Y\to Y'$ and $ j: X \to X'$ such that $i(Y)$ and $j(X)$ are contained in the interior (in the sense of Berkovich \cite[Definition 2.5.7]{Ber90}) of $Y'$ and $X'$, respectively, and such that $\vphi'\circ i=j\circ \vphi$. Usually, we just say that $\vphi'$ is a prolongation and we identify $Y$ and $X$ with their respective images, although we keep in mind that $i$ and $j$ are part of the data. 
 
 When $Y'$ and $X'$ are projective curves, we say that $\vphi'$ is a projectivization of the morphism $\vphi$.
\end{definition}

\begin{theorem}
 There exists a prolongation of the morphism $\vphi$.  
\end{theorem}
\proof
 It is clear that the nontrivial cases are when $Y$ (hence also $X$) are not projective or wide open curves. Assume first that $Y$ and $X$ are strict affinoid curves. We may write $X=X'-\cup_{i=1}^mD'_i$, where $X'$ is a simple projectivization of $X$ and where $D'_i$ are strict open discs in $X'$. Let $x_1,\dots,x_l\in X(k)$ be the branching points of $\vphi$, and for each $x_j$ let $D_j$ be a small enough strict open disc in $X$ which contains $x_j$. Put $X_1:=X-\cup_{j=1}^lD_j$ and let $Y_1=\vphi^{-1}(X_1)$. The restriction $\vphi_1:=\vphi_{|Y_1}$ is a finite \'etale morphism and \cite[Proposition 2.4]{Gar96} (see also Remark below) implies that there exists a finite morphism $\vphi':Y' \to X'$ which extends $\vphi_{1}$. Let $X_2$ be $X'-\cup_{j=1}^lD_j$ and let $Y_2$ be the preimage $\vphi'^{-1}(X_2)$. Then, we can glue  $Y_2$ and $Y$ along $Y_1$ to obtain a (smooth projective) curve $Y''$, and accordingly we can glue morphisms $\vphi_2:=\vphi_{|Y_2}$ and $\vphi$ along $\vphi_{1}$ to obtain a finite morphism $\vphi'':Y''\to X'$ which in fact prolongs  our original morphism $\vphi$. 
 
 If $Y$ and $X$ are strict ft curves that are not smooth or affinoid, then let $X'$ be a simple projectivization of $X$ so that we can write $X=X'-\cup_{i=1}^mB_i-\cup_{j=1}^lD_j$, where $B_i$ are strict closed and  $D_j$ are strict open discs in $X'$. Furthermore, let $X_1$ be a strict connected affinoid domain in $X$ so that the endpoints of discs $D_j$ belong to $X_1$. Moreover, by shrinking $X_1$ if necessary we may assume that the restriction $\vphi_1:=\vphi_{|Y_1}:Y_{X_1}:=\vphi^{-1}(X_1)\to X_1$ is finite and \'etale. Then, again \cite[Proposition 2.4]{Gar96} implies that $\vphi_{1}$ extends to a finite covering of $X_2:=X_1\cup_{j=1}^l D_j$ (note that $X_2$ is an affinoid domain in $X'$) which we denote by $\vphi_2:Y_2\to X_2$. Finally, similarly as before, we can glue $Y_2$ and $Y$ along $Y_1$ and accordingly $\vphi_2$ and $\vphi$ along $\vphi_1$ to obtain an ft curve $Y_3$ and a finite morphism $\vphi_3:Y_3\to X\cup_{j=1}^l D_j$ which is clearly a prolongation  of $\vphi$.  
\qed
\begin{remark}
 In the result \cite[Proposition 2.4]{Gar96} which was used above it is assumed that the base field is discretely valued. However, the same argument from \lc which is based on Krasner's lemma, goes almost verbatim in our situation, so we omit its repetition. 
\end{remark}

\section{Strictly compatible  triangulations of curves}

In order to simplify the study of finite morphisms of affinoid curves, we introduce the notion of strictly compatible triangulations. Later on their existence will be used to deduce the global RH formula from the RH formula for finite morphisms of $k$-affinoid curves which have good reduction. 
\begin{definition}
Let $\vphi:Y\to X$ be a finite morphism of ft curves. We say that (nonempty) triangulations $\cS$ and $\cT$ of $Y$ and $X$, respectively, are $\vphi$-compatible if $\cS=\vphi^{-1}(\cT)$. We say they are strictly $\vphi$-compatible if no two endpoints of an open annulus in $Y-\cS$  are sent to the same point in $\cT$. If $\vphi$ is clear from the context, we may just say that triangulations are compatible or strictly compatible. 
\end{definition}
As one may guess, the compatible triangulations exist.  
\begin{theorem} \label{compatible}
 Let $\vphi:Y\to X$ be a finite morphism of ft curves and let $S\subset Y$ and $T\subset X$ any finite subsets consisting of type two points. Then, there exist strictly $\vphi$-compatible triangulations $\cS$ and $\cT$ of $Y$ and $X$, respectively, such that $S\subset\cS$ and $T\subset\cT$. 
\end{theorem}
\proof
The previous theorem  (or better say its equivalent formulations)  is well known when $Y$ and $X$ are projective curves, and there are several proofs available (\cf \cite{JW,ABBR}, but also \cite{Col03} in the weaker version where one asks for existence of $\vphi$-compatible rather than strictly $\vphi$-compatible triangulations). More precisely, in both articles \cite{JW,ABBR}, authors proved that when $Y$ and $X$ are proper, there exists skeleta $\Gamma_Y$ and $\Gamma_X$ of $Y$ and $X$, respectively, such that $\Gamma_Y=\vphi^{-1}(\Gamma_X)$ (see \cite[Section 3]{JW} or \cite[Theorem A]{ABBR}). Moreover, we may pick skeleta so that they contain any given type two points in $Y$ and $X$, respectively.  Let $\cS'$ and $\cT'$ be triangulations of $Y$ and $X$, respectively, so that $\Gamma_{\cS'}=\Gamma_{Y}$ and $\Gamma_{\cT'}=\Gamma_X$ and $T\subset \cT'$ and $S\subset \cS'$. Then, $\cT:=\vphi(\cS)\cup\cT'\subset \Gamma_X$ is a strict triangulation of $X$ and $\cS:=\vphi^{-1}(\cT)\subset \Gamma_Y$ is a strict triangulation of $Y$ and $\cS$ and $\cT$ are strictly $\vphi$-compatible.

In the case when $Y$ and $X$ are affinoid curves, let $\vphi':Y'\to X'$ be a projectivization of $\vphi$ and let $\cS'$ and $\cT'$ be strictly $\vphi$-compatible triangulations of $Y'$ and $X'$ so that they contain the Shilov boudary of $Y$ and $X$, respectively, and also $S\subset \cS'$ and $T\subset \cT'$. Then, we may take $\cS:=\cS'\cap Y$ and $\cT:=\cT'\cap X$. 

Finally, for $Y$ and $X$ general ft curves, let $Y_0$ and $X_0$ be affinoid subdomains in $Y$ and $X$, respectively, so that $S\subset Y_0$, $T\subset X_0$, $Y-Y_0$ (resp. $X-X_0$) is a finite union of open annuli and such that $\vphi_{|Y_0}:Y_0\to X_0$ is finite. Then, any pair of strictly $\vphi_{Y_0}$-compatible triangulations of $Y_0$ and $X_0$ which contains sets $S$ and $T$, respectively, yields a pair of strictly $\vphi$-compatible triangulations of $Y$ and $X$ which satisfy the claim in the theorem, and we are done.
\qed

We are mainly interested in some consequences of the previous theorem, which we collect in the next corollary.

\begin{corollary} \label{consq strict compa}Keep the notation as in Theorem \ref{compatible}.

{\em (i)} For each open annulus $A\in \cA_{\cT}(X)$, $\vphi^{-1}(A)$ is a disjoint union of open annuli in $\cA_{\cS}(Y)$ and for each $A\in \cA_{\cT}(Y)$, $\vphi(A)\in \cA_{\cT}(X)$ and $\vphi:A\to \vphi(A)$ is a finite morphism. We can furthermore find refinements $\cS\subset\cS'$ and $\cT\subset \cT'$ such that for each $A\in \cA_{\cS'}$, $\vphi_{|A}:A\to \vphi(A)$ is \'etale.

 {\em (ii)} $\Gamma^Y_{\cS}=\vphi^{-1}(\Gamma^X_{\cT})$.

{\em (iii)} Suppose that $\cS$ (hence also $\cT$) has at least two points. For each point $y\in \Gamma^{Y}_{\cS}$, $\vphi$ restricts to a finite morphism $\vphi:C_{\cS,y}\to C_{\cT,\vphi(y)}$. 
 
\end{corollary}
\proof
 {\em (i)} For each open annulus $A\in \cA_{\cS}(Y)$, the image $\vphi(A)$ is contained in some $A'\in \cA_{\cT}(X)$. This is because the triangulations are strictly compatible and $\vphi$ is continuous. But then, one of the connected components of $\vphi^{-1}(A')$ must be $A$, so $\vphi_{|A}:A\to A'$ is finite. Moreover, if $D$ is an open disc in $Y$ with an endpoint $s\in\cS$, $\vphi(D)$ must have an endpoint $\vphi(t)\in\cT$. So, $\vphi(D)$ must be contained in a connected open subset that has an endpoint in $\vphi(s)$. This implies that either $\vphi(D)$ is contained in a residual class that is an open disc attached to $\vphi(s)$, either $\vphi(D)$ is contained in some open annulus $A\in \cA_{\cT}(X)$ such that $A$ has an endpoint in $\vphi(s)$. But again because of strict stability of the triangulations, $D$ is one of the connected components of $\vphi^{-1}(\vphi(D))$ so $\vphi(D)$ is either an open disc attached to $s$ either an open annulus. But, the later is not possible as $\vphi(D\cup\{s\})$ is closed, so $\vphi(D)$ is an open disc attached to $\vphi(s)$. It follows that connected components of $\vphi^{-1}(A)$, for $A\in \cA_{\cT}(X)$ are open annuli in $\cA_{\cS}(Y)$. 
 
 If for some $A\in \cA_{\cS}(Y)$, the restriction $\vphi_{|A}:A\to \vphi(A)$ is not \'etale, let $y_1,\dots,y_m \in A(k)$ be critical points of $\vphi_{|A}$, let $D_{A,1},\dots,D_{A,n(A)}$ be the maximal open discs in $A$ that contain them, and let $\eta_{A,1},\dots,\eta_{A,n(A)}$ be the endpoints of $D_{A,1},\dots,D_{A,n(A)}$, respectively. Then, $\eta_{A,1},\dots,\eta_{A,n(A)}$ lie on the skeleton $\Gamma^A$, and the refinement that we are looking for is $\cT':=\cT\cup_{A\in \cA_{\cS}(Y)}\{\vphi(\eta_{A,1},\dots,\vphi(\eta_{A,n(A)})\}$ and $\cS':=\vphi^{-1}(\cT')$.
 
 {\em (ii)} is a direct consequence of {\em (i)} and Lemma \ref{skeleta annuli}, while for {\em (iii)} we note that because of {\em (ii)} each open disc in $C_{\cS,y}$ which is attached to $y$ is sent to an open disc in $C_{\cT,\vphi(y)}$. In particular, $C_{\cS,y}$ is a connected component of $\vphi^{-1}(C_{\cT,\vphi(y)})$ and the result follows.
\qed

\section{De Rham cohomology of ft curves}

 To define Euler-Poincar\'e characteristic of an ft curve we will use the de Rham cohomology with respect to its overconvergent structure sheaf, and in this section we recall some of its basic properties. We refer to \cite{GK00,GK04} and in particular to \cite{VDPcoh} for the main notions and results concerning de Rham cohomology on affinoid spaces. The main result is Theorem \ref{dR groups} which gives us dimensions of the de Rham cohomology groups of ft curves in dependence on  the dimensions of the de Rham cohomology groups of their simple projectivizations and  the number of "missing" discs.

\subsection{De Rham cohomology of strict ft curves}
To calculate the de Rham cohomology groups on strict ft curves, we equip them with overconvergent structure sheaf, as in \cite{GK04}. More precisely, if $X$ is an ft curve and $X'$ a projectivization of $X$, we define $O_{X'}(X):=\limind_{X\subset V}O_{X'}(V)$, where $V$ goes through a set of open neighborhoods of $X$ in $X'$. The $k$-algebra $O_{X'}(X)$ in fact does not depend on the chosen projectivization of $X$, as can be shown for example using \cite[Corollary 1]{Bos81}. We put $O(X^\dg)$ for $O_{X'}(X)$ if we don't need to emphasize the projectivization $X'$. Similarly as before, we let $d:O(X^\dg)\to \Omega^1_{X^\dg}$ to be the universal $k$-linear derivation of $O(X^\dg)$ into finite $O(X^\dg)$-modules (\cf \cite[Section 4.1]{GK00}).
\begin{definition}
 We define $\H^i_{dR}(X)$, $i=0,1,2$ to be the $i^{th}$ hypercohomology group of the De Rham complex $0\to O(X^\dg)\to \Omega^1_{X^\dg}$. 
\end{definition}
\begin{remark}\label{coh remarks}
 $a)$ If $X$ is projective, then the De Rham cohomology groups are isomorphic to algebraic De Rham cohomology groups because of $p$-adic GAGA.\\
 $b)$ If $X$ is a strict ft curve that is not projective, then it is a quasi-Stein space, \cf \cite[Definition 2.3]{Kie67}. Recall that the latter (roughly) means that there exists a sequence of affinoid domains $U_i$, $i\in \N$ in $X$, whose union covers $X$  and such that: 1) $U_i\subseteq U_{i+1}$ and 2) The affinoid algebra $\cO(U_{i+1})$ is dense in the affinoid algebra $\cO(U_i)$ with respect to the spectral norm. To show our claim, we may assume that $X$ is not compact (otherwise we may take $U_i=X$). Let $\cT$ be a nonempty triangulation of $X$ and let $A_1,\dots,A_n$ be the open annuli in $\cA_{\cT}$ which are not relatively compact in $X$. For each $j=1,\dots,n$, let $A_{j,i}$, $i\in \N$ be a sequence of strict open annuli in $A_j$ such that $A_{j,i+1}\subset A_{j,i}$, $\bigcap_{i\in\N}A_{j,i}=\emptyset$ and for every $i\in \N$, $A_{j,i}$ is not relatively compact in $X$. Then, we may take $U_i:=X-\Big(\bigcup_{j=1}^n A_{j,i}\Big)$, $i\in \N$. The condition 1) is satisfied by the construction of $U_i$ while the condition 2) is a consequence of $p$-adic Runge's theorem, a result that we will encounter later in Theorem \ref{runge} (see Remark \ref{coh remarks continued}).
 
 As a consequence, the higher cohomology groups for coherent $O_X$-modules vanish (\cf \cite[Proposition 3.1]{GK00}) and we can calculate the de Rham cohomology groups directly as the cohomology groups of the global sections of the de Rham complex. This fact will be used in what follows. 
 \end{remark}

\begin{example} \label{exa dR}
 $(a)$ Let $X$ be a strict open disc, and let $T:X\iso D(0,1^-)$ be a coordinate. Then $O(X^\dg)$ can be identified with the ring of functions 
 $$
\Big \{f(T)=\sum_{i\geq 0} a_i T^i\in k[[T]],\quad \lim_{i\to\infty}|a_i|\rho^i=0,\quad \forall \rho\in(0,1)\Big\},
 $$
or with the power series $f(T)\in k[[T]]$ which have the radius of convergence greater than or equal to 1, that is we have $O(X^\dg)=O(X)$. Let us denote the latter ring temporarily by $D(T)$. Then, $\Omega^1_{X^\dg}$ can be identified with $D(T)dT$ and in this case $d=\frac{d}{dT}$. It is a direct computation to show that $\H^0_{dR}(X)\simeq k$ and $\H^1_{dR}(X)\simeq 0$ as for any $f(T)\in D(T)$ and $g(T)\in k[[T]]$ with $g'(T)=f(T)$ one can check that $g(T)\in D(T)$. 

$(b)$ Similarly, if $X$ is a strict open annulus, and $T:X\iso A(0;r,1)$ a coordinate, then $O(X^\dg)$ can be identified with the ring of series 
$$
\Big\{f(T)\in\sum_{i\in \Z}a_i T^i,\quad \lim_{|i|\to\infty} |a_i|\rho^i=0,\quad \forall \rho\in(r,1)\Big\}.
$$
and if we denote the later by $A(T)$, then $\Omega^1_{X^\dg}$ can be identified with $A(T)dT$ with $d=\frac{d}{dT}$. Once again, a direct computation shows that $\H^0_{dR}(X)\simeq k$ while $\H^1_{dR}(X)\simeq k$, because $\frac{1}{T}$ is not integrable in $A(T)$. We note that the same conclusions hold when $X$ is a strict punctured open disc.
 \end{example}

 The following lemma is implicitly contained in \cite[Section 4]{Col89} (but see also \cite[Section 2]{ColemanRegulators}).

\begin{lemma}\label{coh wide open}
 Let $X$ be a strict wide open curve which is not projective, let $X'$ be its simple projectivization, and let $\delta$ be the number of connected components in $X'-X$. Then $\H^0_{dR}(X)\simeq k$ and $\H^1_{dR}(X)\simeq k^{2g+\delta-1}$, where $g$ is the genus of $X'$, while $\H^i_{dR}(X)=0$, for $i\geq 2$.
\end{lemma}
\proof
 Let $D'_i$, $i=1,\dots,\delta$ be the connected components of $X'-X$, and for each $i$ let $a_i$ be a rational point in $D'_i$. Let $X'':=X'-\{a_1,\dots,a_\delta\}$, so that $X''$ is the analytification of a smooth affine $k$-algebraic curve denoted by $\cX''$. For each $i$ let $A_i$ be a strict open annulus (in fact an open disc punctured in the point $a_i$) in $X''$ with the following properties: $A_i$ intersects $X$ and the intersection is a strict open annulus denoted by $A''_i$, $A_i\cap A_j=\emptyset$ for $i\neq j$ and finally $X$ and the union of $A_i$ cover $X''$ (such $A_i$ exist because $X$ is a wide open curve). Then, we have the Mayer-Vietoris sequence for the hypercohomology:
 \begin{align}
 0&\to \H^0_{dR}(X'')\to \H^0_{dR}(X)\oplus\bigoplus_{i=1}^\delta \H^0_{dR}(A_i)\to \oplus_{i=1}^\delta\H^0_{dR}(A''_i)\to \nonumber\\
&\to  \H^1_{dR}(X'')\to \H^1_{dR}(X)\oplus\bigoplus_{i=1}^\delta \H^1_{dR}(A_i)\to \oplus_{i=1}^\delta\H^1_{dR}(A''_i)\to \label{MV seq}\\
&\to \H^2_{dR}(X'')\to \H^2_{dR}(X)\to 0.\nonumber
 \end{align}
 At this point we recall that $\H^i_{dR}(X'')\simeq \cH^i_{dR}(\cX'')$ where the later is the algebraic de Rham cohomology of $\cX''$ and this implies for $i\geq 2$, $\H^{i}_{dR}(X'')=\H^i_{dR}(X)=0$. As $\H^0_{dR}(X)\simeq k$ naturally, and having in mind the Example \ref{exa dR} $(a)$ and $(b)$,  for dimension reasons the sequence \eqref{MV seq} decomposes in two short exact sequences. The first, which corresponds to the first row of \eqref{MV seq}, becomes $0\to k\to k\oplus k^\delta\to k^\delta\to 0$, while the second one is $0\to k^{2g+\delta-1}\to \H^1_{dR}(X)\oplus k^{\delta}\to k^\delta\to0$ which implies $\H^1_{dR}(X)\simeq k^{2g+\delta-1}$. 
\qed
\begin{lemma}\label{cov}
 Let $U$ be a strict wide open curve and let $\{U_i\}_{i\in I}$ be a finite covering of $U$ by strict wide open curves such that $U_i\cap U_j\cap U_l=\emptyset$ whenever $i\neq j$, $i\neq l$ and $j\neq l$. For a subset $J\subset I$ let us denote by $U_J:=\bigcap_{j\in J}U_j$ and by $|J|$ the cardinality of $J$. Then, we have an exact sequence 
 \begin{align*}
 0\to \H^0_{dR}(U) &\to \bigoplus_{J\subset I,|J|=1}\H^0_{dR}(U_J)\to\bigoplus_{J\subset I} \H^0_{dR}(U_J)\to \\
 &\to \H^1_{dR}(U)\to \bigoplus_{J\subset I,|J|=1}\H^1_{dR}(U_J)\to\bigoplus_{J\subset I} \H^1_{dR}(U_J)\to 0.
 \end{align*}
\end{lemma}
\proof
 Let us for the moment denote by $\cF^*$ the overconvergent De Rham complex on $U$, and for $J\subset I$, let $j_J:U_J\to U$ be the inclusion of $U_J$ into $U$ and let $\cF^*_J=j_{J*}(\cF^*_{|U_J})$. First we observe that the sequence of complexes $0\to \cF^*\to \bigoplus_{|J|=1}\cF^*_J\to \bigoplus_{|J|=2}\cF^*_J\to 0$ is exact, as for each $i\geq0$, $0\to \cF^i\to \bigoplus_{|J|=1}\cF^i_J\to\bigoplus_{|J|=2}\cF^i_{J}\to 0$ is a \v Cech resolution of $\cF^i$ with respect to the covering $\{U_i\}_{i\in I}$. Then, such an exact sequence of complexes induces a long exact sequence of the form 
$$
0\to\H^0(U,\cF^*)\to \bigoplus_{J\subset I, |J|=1} \H^0(U, \cF_J^*)\to 
\bigoplus_{J\subset I, |J|=2} \H^0(U,\cF_J^*) \to 
\H^1(U,\cF^*)\to\cdots,
$$
where the groups involved are just hypercohomology groups for the corresponding complexes. Finally, we note that $\H^l(U,\cF^*)=\H^l_{dR}(U)$ and $\H^l(U,\cF^*_J)=\H^l(U_J,\cF^*_{|U_J})=\H^l_{dR}(U_J)$ because $U_J$ and $U$ are quasi-Stein. 
\qed

Now we are ready to calculate the De Rham cohomology groups of general ft curves. The argument is also based on \cite{GK13'}.
\begin{theorem}\label{dR groups}
 Let $X$ be an ft curve, let $X'$ be its simple projectivization, and let $m$ be the number of connected components in $X'-X$. Then, if $m=0$ we have $H^0_{dR}(X)\simeq k$, $H^1_{dR}(X)\simeq k^{2g}$, $H^2_{dR}(X)\simeq k$ and $H^i_{dR}(X)\simeq 0$ for $i\geq 3$. Otherwise, $\H^0_{dR}(X)\simeq k$, $\H^1_{dR}(X)\simeq k^{2g+m-1}$ and $\H^i_{dR}(X)=0$ for $i\geq 2$.  
\end{theorem}
\proof
 If $m=0$ then the theorem follows from Remark \ref{coh remarks} and corresponding isomorphisms for algebraic De Rham cohomology groups. Suppose $m>0$. We may also assume that $X$ is not a wide open curve as this case was done in Lemma \ref{coh wide open}. Then, there exists a strictly decreasing sequence $(U_n)_{n\geq 0}$ of strict wide open curves in $X'$ such that: 1) $\cap_n U_n=X$ and 2) for each $n\geq 0$, $U_n-X$ is a disjoint union of strict open annuli (note that the number of these open annuli is the same for each $n$). Let $U_0-X=\uplus_{t=1}^l A_t$ and for each $n\geq 1$ and $t=1,\dots,l$ put $A_{n,t}:=A_t\cap U_n$. Note that for each $n\geq 1$, the finite family $\{U_n, A_{1},\dots A_l \}$ is an open covering of $U_0$ which satisfies the conditions of Lemma \ref{cov}, so that we have an exact sequence 
  \begin{align*}
 0\to \H^0_{dR}(U_0) &\to \H^0_{dR}(U_n)\oplus\bigoplus_{t=1}^l\H^0_{dR}(A_t)\to\bigoplus_{t=1}^l \H^0_{dR}(A_{n,t})\to \\
 &\to\H^1_{dR}(U_0)\to \H^1_{dR}(U_n)\oplus\bigoplus_{t=1}^l\H^1_{dR}(A_t)\to\bigoplus_{t=1}^l \H^1_{dR}(A_{n,t})\to 0.
 \end{align*}
 Taking the direct limit $\limind_{n}$ and having in mind that direct limit commutes with taking the cohomology groups (so that we have $\H^i_{dR}(X)=\limind_n\H^i_{dR}(U_n)$), we obtain
  \begin{align*}
 0\to \H^0_{dR}(U_0) &\to \H^0_{dR}(X)\oplus\bigoplus_{t=1}^l\H^0_{dR}(A_t)\to\bigoplus_{t=1}^l \limind_{i}\H^0_{dR}(A_{n,t})\to \\
 &\to\H^1_{dR}(U_0)\to \H^1_{dR}(X)\oplus\bigoplus_{t=1}^l\H^1_{dR}(A_i)\to\bigoplus_{t=1}^l \limind_n\H^1_{dR}(A_{n,t})\to 0.
 \end{align*}
 Finally, we note that for each $t=1,\dots,l$, $\dim_k \limind_n \H^0_{dR}(A_{n,t})=\dim_k \limind_n \H^1_{dR}(A_{n,t})=1$ and the result follows by comparing the dimensions and using Lemma \ref{coh wide open}. 
\qed
\begin{remark}
 One may notice that $\limind_n \H^*_{dR}(A_{n,t})$ corresponds to the De Rham cohomology groups of Robba annuli $\limind_nA_{n,t}$.
\end{remark}

\subsection{Euler-Poincar\'e characteristic of an ft curve} 
\begin{definition}
Let $X$ be an ft curve. The number $\chi(X):=\sum_{i=0}^2(-1)^i\dim_k\H^i_{dR}(X)$ is called the Euler-Poincar\'e characteristic of $X$. 
\end{definition}
With similar arguments as in Theorem \ref{dR groups} and the fact that Euler-Poincar\'e characteristic of an open annulus is 0 (see Example \ref{exa dR} $(b)$) one can prove the additivity of Euler-Poincar\'e characteristic, in the sense of the next corollary. We state it without proof.
\begin{corollary}\label{EP add}
 Let $X$ be an ft curve, and let $\cT$ be a triangulation of $X$ which contains at least two points. Then, $\chi(X)=\sum_{x\in \cT}\chi(C_{\cT,x})$. 
\end{corollary}

\section{Riemann-Hurwitz formula}

\subsection{Finite \'etale morphisms of open annuli}

\subsubsection{}
Let $\vphi:A_1\to A_2$ be a finite \'etale morphism 
of open annuli of degree $d$. Let $S:A_2\iso A(0;\rho^d,1)$ 
(resp. $T:A_1\iso A(0;\rho,1)$) be a coordinate on $A_1$ 
(resp. $A_2$). Then $\vphi$ can be represented as $S=\vphi_{\#}(T)=a_dT^du(T)$, where 
$|u(T)-1|_{\rho_0}<1$ for every $\rho_0\in (\rho,1)$ and with $|a_d|=1$. By choosing a different coordinate $T$, we may achieve that $a_d=1$ so we will assume, unless otherwise stated that $a_d=1$. Since $\vphi$ is \'etale, the derivative of $\vphi$ is an 
invertible function on $A_1$, hence has the following coordinate representation 
$\frac{dS}{dT}=\vphi_{\#}'(T)=\eps T^\sigma v(T)$, where again 
$|v(T)-1|_{(\rho_0)}<1$, $\rho_0\in (\rho,1)$ and where $\eps\in \kc$. We put 
\beq\label{nu term}
\nu=\sigma-d+1.
\eeq
If we want to emphasize the dependence on the morphism $\vphi$, we will write $\sigma(\vphi)$, $d(\vphi)$, $\nu(\vphi)$,  $\eps(\vphi)$, and so on. 

\begin{remark}
 1) To the best of author's knowledge, the invariant $\nu$ first time appears in the paper \cite{Lut93} by W. L\"utkebohmert where it is related to certain discriminant associated to coverings of closed annuli. To be more precise, let us keep the notation as above and let $r\in (\rho,1)\cap |k|$.  We note that $\vphi$ restricts to a finite \'etale morphism $\vphi_r:A[0;r,r]\to A[0;r^d,r^d]$ and in this case one can study the discriminant of the extension $\vphi_r$ (\lc Definition 1.3) as well as the absolute value of the discriminant which is in \lc denoted by $\mathfrak{d}(r^d)$ (in fact, $\mathfrak{d}$ is a function of the radius of the target annulus). Then Lemma 1.7 in \lc shows that   
 $$
 \mathfrak{d}(r^d)=\begin{cases}
                    |d|^d,\quad \nu=0;\\
                    |\eps|^d\cdot r^{\nu d},\quad\text{ otherwise.}
                   \end{cases}
 $$
 One can extend the function $\mathfrak{d}$ to a continuous function on the whole segment $(\rho^d,1)$ in an obvious way (or, one can extend the definition of $\mathfrak{d}$ to points $r\notin(\rho,1)\cap|k|$ by a suitable extension of scalars!) and in this case $\nu$ is precisely the logarithmic slope of the function $\mathfrak{d}$.

 2) A different paradigm was adopted in \cite{CTT14} where authors introduced the different function $\delta_\vphi$. Roughly speaking, if $\vphi:Y\to X$ is a finite morphism of ft curves and $y\in Y$, then $\delta_\vphi(y)$ is the absolute value of the different ideal of the extension $\sH(y)/\sH(\vphi(y))$. If we are in the situation above and $r\in (\rho,1)$, then one shows that $\delta_\vphi(\eta_{0,r})=|\eps|\cdot r^\nu$, so in this context our invariant $\nu$ is simply the logarithmic slope of the different function along the skeleton of $A(0;\rho,1)$. One may note that the different function is exactly $\deg(\vphi)$-root of the discriminant function from above, as one would usually expect.
 
 3) If $r\in(\rho,1)\cap|k|$, then the residue field extension $\wtilde{\sH(\eta_{0,r})}/\wtilde{\sH(\eta_{o,r^d})}$ (here $\eta_{0,r^d}=\vphi(\eta_{0,r})$) is separable if and only if $|\eps|\cdot r^\nu=1$ (\cite[Lemma 4.2.2.]{CTT14}). In this case, which happens for example when $\char(\kt)=0$, we have $\nu=0$.
 \end{remark}

\subsubsection{} 
We collect some of the properties of finite \'etale morphisms of open annuli that will be used later on.

\begin{lemma} \label{rel annuli}
 Let $\vphi:A(0;\rho,1)\to A(0;\rho^{d(\vphi)},1)$ and 
$\psi:A(0;\rho^{d(\vphi)},1)\to A(0;\rho^{d(\vphi) d(\psi)},1)$ be finite 
\'etale morphisms 
of  open annuli of degree $d(\vphi)$ and $d(\psi)$, respectively. Then
\begin{itemize}
\item[(i)] $|\epsilon(\vphi)|\geq|d(\vphi)|$ (here $|\cdot|$ is the norm of the field $k$).
 \item[(ii)] Let $S_1=r^{d(\vphi)}/S$ and $T_1=r/T$, where $r\in k$, $|r|=\rho$ be "inverted" coordinates on $A(0;\rho,1)$ and $A(0;\rho^{d(\vphi)},1)$, respectively. Let $\sigma_1(\vphi)$ (and similarly $\eps_1(\vphi)$, $\nu_1(\vphi)$) be the order of the derivative of the function $\vphi$ expressed in new coordinates $S_1$ and $T_1$.  Then, 
$$
 \sigma_1(\vphi)=-\sigma(\vphi)+2d(\vphi)-2,\quad
 \eps_1(\vphi)=r^{\nu(\vphi)}\cdot \eps(\vphi),\quad \text{ and } \quad
 \nu_1(\vphi)=-\nu(\vphi).
$$

\item[(iii)] We have $\sigma(\psi\circ\vphi)=d(\vphi)\sigma(\psi)+\sigma(\vphi)$ and $\nu(\psi\circ\vphi)=d(\vphi)\nu(\psi)+\nu(\vphi)$.
\end{itemize}
\end{lemma}

\proof
$(i)$ We can write $\vphi_{\#}(T)=\sum_{l\in \Z}a_l T^l$, so 
$\frac{d}{dT}\vphi_{\#}(T)=\sum_{l\in \Z}la_lT^{l-1}$, which implies that 
$\epsilon(\vphi)=a_{\sigma_l+1}(\sigma_{l}+1)$. Then 
$|\epsilon(\vphi)|=|a_{\sigma_l+1}(\sigma_{l}+1)|\geq |d(\vphi)a_{d(\vphi)}|=|d(\vphi)|$. 
For $(ii)$, we have
$$
S_1=\frac{r^{d(\vphi)}}{\vphi(\frac{r}{T_1})}=T_1^{d(\vphi)} (u(\frac{r}{T_1}))^{-1}=T_1^{d(\vphi)} 
u_1(T_1)
$$
and
\begin{align*}
\frac{dS_1}{dT_1}&=r^{d(\vphi)}\cdot\frac{1}{\big(\vphi(\frac{r}{T_1})\big)^2}
\cdot\vphi_{\#}'(\frac{r}{T_1})\cdot\frac{r}{T_1^2}
=r^{d(\vphi)}\cdot\frac{T_1^{2d(\vphi)}}{r^{2d(\vphi)}}\cdot(u_1(T))^2\cdot\eps(\vphi)\frac{r^{\sigma(\vphi)}}{T_1^{
{\sigma}(\vphi)}}\cdot v(\frac{r}{T_1})\cdot\frac{r}{T_1^2}\\
&=r^{\nu(\vphi)}\cdot\eps(\vphi)\cdot T_1^{-\sigma(\vphi)+2d(\vphi)-2}\cdot v_1(T_1)
=\eps_1(\vphi)\cdot T_1^{\sigma(\vphi)}\cdot v_1(T_1).
\end{align*}
We just note here that $u_1(T_1)=u(r/T_1)$ and 
$v_1(T_1)=\big(u_1(T_1)\big)^2v(r/T_1)$ are units so formulae follow.

$(iii)$ If we introduce coordinates $U$, $S$ and 
$T$ on $A(\rho^{d(\vphi) d(\psi}),1)$, $A(\rho^{d(\vphi}),1)$ and $A(\rho,1)$, 
respectively, we may write $U=\psi_{\#}(S)=S^{d(\psi)}h_1(S)$, and 
$S=\vphi_{\#}(T)=T^{d(\vphi)}h_2(T)$, where $h_1$ and $h_2$ are units in their 
respective rings, and $|h_2|_{\rho_0}<1$ for $\rho_0\in (\rho,1)$ and $|h_1|_{\rho_0}<1$ for $\rho_0\in (\rho^{d(\vphi)},1)$ . Let us write $\frac{dU}{dS}=\epsilon(\psi)S^{\sigma(\psi)}g_1(S)$ and $\frac{dS}{dT}=\epsilon(\vphi)T^{\sigma(\vphi)}g_2(T)$ with usual assumptions on functions $g_1$ and $g_2$. Then, it 
is a straightforward computation using the chain rule:
\begin{align*}
\frac{dU}{dT}&=\frac{d}{dT}(\psi_{\#}(\vphi_{\#}(T)))=\frac{d\psi_{\#}}{dS}(\vphi_{\#}(T))\frac{d\vphi_{\#}}{dT}(T)\\
&=\epsilon(\psi)(\vphi_{\#}(T))^{\sigma(\psi)}g_1(\vphi_{\#}(T))\epsilon(\vphi)T^{\sigma(\vphi)}g_2(T)\\
&=\epsilon(\vphi)\epsilon(\psi) T^{d(\vphi)\sigma(\psi)+\sigma(\vphi)}g_1(\vphi_{\#}(T))g_2(T),
\end{align*}
which implies
$\sigma(\psi\circ\vphi)=d(\vphi)\sigma(\psi) +\sigma(\vphi)$. Then,
$\nu(\psi\circ\vphi) =\sigma(\psi\circ\vphi)-d(\vphi) 
d(\psi)+1=d(\vphi)(\sigma(\psi)-d(\psi)+1)+\sigma(\vphi)-d(\vphi)+1=d(\psi)\nu(
\vphi) +\nu(\psi)$.
\qed
\begin{remark}
Although in the beginning of this section we made an implicit dependence of terms $\sigma$ and $\nu$ on the coordinates, the previous lemma and in particular their behavior with respect to composition of morphisms, shows that a change of coordinates does not affect them. 
\end{remark}
The order of derivative $\sigma$ has a very close relation with ramification, as is shown in the next lemma. A generalization is given in Theorem \ref{sigma general}.
\begin{lemma} \label{sigma ram}
 Suppose that $\vphi:A_1\to A_2$ extends to a finite map of the whole open disc 
$D(0,1^-)$ to itself with ramified points $x_1,\dots,x_s\in D(0,1^-)(k)$ and ramification 
indexes 
$e_{x_1},\dots,e_{x_s}$, respectively. Then $\sigma=\sum_{1\leq 
i\leq s}(e_{x_i}-1)$.
\end{lemma}
\proof
 The map  $\vphi$ has a coordinate representation as a power 
series $S=\vphi_{\#}(T)=\sum_{i\geq0}a_i T^i$, and the derivative of $\vphi$ 
is again a power series that we can factor as $\frac{dS}{dT}=P(T)g(T)$, where 
$P(T)=(T-x_1)^{e_{x_1}-1}\cdots(T-x_s)^{e_{x_s}-1}$, 
while $g(T)$ is invertible on $D(0,1^-)$. As 
$\rho\in (0,1)$ approaches $1$, the theory of valuation polygons tells us that 
the the logarithmic derivative $\text{dlog}^-|\frac{dS}{dT}|_{\rho}(1)$ is, on 
one side equal to the number of zeros in $D(0,1^-)$ of $\frac{dS}{dT}$ counted 
with multiplicities, so exactly $\sum_{1\leq i\leq s}(e_{x_i}-1)$, while on the 
other side it is equal to $\sigma$.
\qed

By finding suitable coordinate representation of the morphism $\vphi$, one can without many difficulties prove the following:
\begin{corollary}
 Suppose that $\vphi$ is ramified at exactly one point in $D(0,1^-)$. Then $\sigma\leq d-1$. 
\end{corollary}

\subsubsection{ The norm of the operator $\frac{d}{dT}$}
Let $T:A\iso A(0;\rho,1)$ be a coordinate on an open annulus $A$. Then, 
every function $f$ on $A$ can be seen via 
$T$ as a function $f(T)=\sum_{i\in \Z}a_i T^i$, where coefficients $a_i\in k$ 
satisfy the condition: for each 
$\rho_0\in(\rho,1)$, $\lim_{|i|_{\infty}\to \infty}|a_i|\rho_0^i=0$, where 
$|\cdot|_{\infty}$ is the usual archimedean absolute value. 

On the other side, the derivative $\frac{d}{dT}f(T)=\sum_{i\in \Z}a_i i 
T^{i-1}$ can also be seen as a function on $A(0;\rho,1)$ because for all 
$i\in\Z$, $|i|\leq 1$ and therefore for each $\rho_0\in(\rho,1)$,  
$\lim_{|i|_{\infty}\to \infty}|a_i||i|\rho_0^{i-1}= 0$. In this way we can 
see $\frac{d}{dT}$ as an operator acting on the space of fuctions on 
$A(0;\rho,1)$.

\begin{lemma}\label{derivative norm}
Let $\rho_0\in(\rho,1)$ and let $|\frac{d}{dT}|_{\rho_0}$ be the operator norm 
of the operator $\frac{d}{dT}$ seen as acting on the space of functions on 
$A(0;\rho,1)$ equipped with the norm $|\cdot|_{\rho_0}$. Then, 
$|\frac{d}{dT}|_{\rho_0}=\rho_0^{-1}$.
\end{lemma}
\proof
 The proof is straightforward: if $f(T)=\sum_{i\in\Z}a_i T^i$ is an analytic 
function on $A(0;\rho,1)$, then 
$\frac{d}{dT}f(T)=\sum_{i\in\Z}a_i i T^{i-1}$. Furthermore, 
$$
|T\frac{d}{dT}f(T)|_{\rho_0}=\max_{i\in\Z}|a_i||i|\rho_0^i\leq 
\max_{i\in\Z}|a_i|\rho_0^i=|f(T)|_{\rho_0}
$$
which implies 
$|\frac{d}{dT}f(T)|_{\rho_0}\leq \rho_0^{-1}|f(T)|_{\rho_0}$. Since 
$|\frac{d}{dT}(T)|_{\rho_0}=\rho_0^{-1}|T|_{\rho_0}$, the proof follows.
\qed

\subsubsection{The $\sigma(\vphi,\vt)$ and $\nu(\vphi,\vt)$ terms}\label{sigma nu}

Given a morphism $\vphi:Y\to X$ of ft curves, $y\in Y$ a type 2 point, and $\cS$ and $\cT$ strictly $\vphi$-compatible triangulations so that $y\in \cS$, it follows from Corollary \ref{consq strict compa} that there exists an affinoid domain $C_y$ in $Y$ such that $\vphi_y=\vphi_{|C_y}:C_y\to \vphi(C_y)$ is a finite morphism and such that $C_y$ has good reduction and  $y$ is its Shilov point. The degree of the morphism $\vphi_{|C_y}$ is independent of the $C_y$ chosen, and is denoted by $\deg(\vphi,y)$. Let $x=\vphi(y)$ and let $C_x:=\vphi(C_y)$. For each element $\vt\in T_yY$ and the corresponding component $U_{\vt}$ in $Y-\cS$, $\vphi$ induces a finite morphism $\vphi_{U_{\vt}}:U_{\vt}\to \vphi(U_{\vt})$. 

If $U_{\vt}$ is an open disc, then so is $\vphi(U_{\vt})$ and both can be normalized with a pair of coordinates $T:U_\vt\iso D(0,1^-)$ and $S:\vphi(U_{\vt})\iso D(0,1^-)$ and $\vphi_{|U_{\vt}}$ expressed with coordinates $S$ and $T$ is a power series $S=\vphi_\#(T)$. The order of the function and of its derivative near the boundary point are denoted by $d(\vphi,\vt)$ and $\sigma(\vphi,\vt)$, respectively. The terms $d(\vphi,\vt)$ and $\sigma(\vphi,\vt)$ are independent of the pair of normalizing coordinates $S$ and $T$. 

Similarly, if $U_{\vt}$ is an open annulus, then so is $\vphi(U_{\vt})$ and we can normalize them with a pair of coordinates $T:U_{\vt}\iso A(0;r,1)$ and $S:\vphi(U_{\vt})\iso A(0;r^\alpha,1)$. The morphism $\vphi_{|U_{\vt}}$ expressed in coordinate form is a series whose order is denoted by $d(\vphi,\vt)$ and the order of its derivative near the boundary point corresponding to $y$ is denoted by $\sigma(\vphi,\vt)$. Both terms are independent of the pair of normalizing coordinates. 

Finally, if $y$ is in the Shilov boundary of $Y$, and if $\vt\in T_yY^\dg$, then after prolonging our morphism with $\vphi':Y'\to X'$, where $Y'$ and $X'$ are some wide open curves, $\vphi'$ induces a finite (\'etale) morphism of some small enough strict open annuli $U_{\vt}$ and $\vphi'(U_{\vt})$. In this case, similarly as before, we define terms $\sigma(\vphi,\vt)$ and $d(\vphi,\vt)$. However, we must show the independence of the chosen prolongation.
\begin{lemma}
 Let $\vphi_i:Y_i\to X_i$ be two prolongations of the morphism $\vphi:Y\to X$ of quasi-smooth affinoid curves. Let $\vt\in TX$. Then, $d(\vphi_1,\vt)=d(\vphi_2,\vt)$ and $\sigma(\vphi_1,\vt)=\sigma(\vphi_1,\vt)$.
\end{lemma}
\proof
 The claim is local around points of the Shilov boundary of $Y$, so let $\eta$ be a one such point and put $\xi=\vphi(\eta)$. Finite morphism $\vphi_i$ induces a finite \'etale extension of complete residue fields $\sH(\xi)\hookrightarrow\sH(\eta)$ and the extension is independent whether we consider $\eta \in Y_1$ or $\eta\in Y_2$. Theorem \cite[Theorem 3.4.1]{BerCoh} implies that there is an isomorphism $\alpha$ of some open neighborhoods $\cU_1$ and $\cU_2$ of $\eta$ in $Y_1$ and $Y_2$ such that $\vphi_{1|\cU_1}$ and $\vphi_{2|\cU_2}$ are finite morphism and such that $\vphi_{2|\cU_2}\circ\alpha=\vphi_{1|\cU_1}$. The fact that $\alpha$ is an isomorphism which implies that its restriction to any residue class at $\eta$ is an isomorphism as well, together with Lemma \ref{rel annuli} $(iii)$ and remark that follows it yields the claim.
\qed

\begin{definition}\label{nu term gen}
If $\vphi:Y\to X$ is a finite morphism of ft curves, $y\in Y$ a type two point and $\vt\in T_yY^\dg$, we define  $\nu(\vphi,\vt):=\sigma(\vphi,\vt)-d(\vphi,\vt)+1$, and where $\sigma(\vphi,\vt)$ and $d(\vphi,\vt)$ are defined above.
\end{definition}
\subsubsection{Reduction of morphisms} \label{red morph}
On the other side, the reduction of $\vphi_y$ induces a finite morphism of smooth $\kt$-algebraic curves $\wtilde{\vphi}_y:\wtilde{C}_y\to \wtilde{C}_x$. If $\vt\in \wtilde{C}_y$ is a smooth point then the multiplicity of $\wtilde{\vphi}$ at $\vt$ is equal to the degree of the morphism $\vphi_{U_\vt}$, constructed above and we denote this number by $\deg(\vphi,\vt)$.  For a more detailed study and proofs of the previous statements we refer to \cite{Duc-book} or \cite[Section "Coordinates and tangent space"]{Bal-Ked}. 

A direct consequence of the previous remark is the following result that will be used in the proof of Riemann-Hurwitz formula.
\begin{lemma}\label{pullback divisor}
 Let $\vphi_y:C_y\to C_x$ be as above, and let $D=\sum_i a_i\vv_i$ be a divisor on $\wtilde{C}_x$. Then, 
 $$
 \wtilde{\vphi}_y^*(D)=\sum_{i}\sum_{\substack{\vt\in \wtilde{C}_y\\ \wtilde{\vphi}_y(\vt)=\vv_i}}\deg(\vphi_y,\vt)\vt \text. 
  $$
\end{lemma}

\subsection{$p$-adic Runge's theorem}
One of the main ideas behind our proof of RH formula is that certain morphisms/functions on ft curves can be well approximated by rational functions on the corresponding simple projectivizations, \ie the $p$-adic Runge's theorem.
\begin{theorem}\label{runge}
Let $Y$ be a strict, quasi-smooth $k$-affinoid curve, and let $Y'$ be its simple projectivization. For each connected component $D_i$, $i=1,\dots,m$ in $Y'-  Y$ let $y_i\in D_i(k)$. Then, the ring of rational functions on $Y'$ and with possible poles only at points $y_1,\dots,y_m$ is dense in the ring on holomorphic functions on $Y$ with respect to the spectral norm on $Y$.
\end{theorem}
\proof
Let $\cY'$ be a smooth projective curve such that the analytification $\cY'^{an}$ of $\cY'$ is isomorphic to $Y'$. Let us still denote by $y_i$ the corresponding points in $\cY'$ and let $\cY:=\cY'-\{y_1,\dots,y_m\}$. If we can find an affine $k^\circ$-scheme of finite type $\cU$ such that: 1) $\cU\otimes k=\cY$ and 2) $\widehat{\cU}_\eta$- the generic fiber of the formal completion of $\cU$- is isomorphic to $Y$, then we are done. The existence of such a scheme $\cU$ can be proved along the same lines as in \cite[Proposition 3.5.1]{Ray94} with slight changes, so we just sketch the proof. Namely, let us start with any (semi-stable) $k^\circ$-model $\cU'$ of $\cY'$ (see \cite[Section 1]{Bal10}), such that if $spe:\widehat{\cU'}_\eta\to \widehat{\cU'}_s=\cU'_s$ is the specialization map, then the Shilov boundary of $Y$ is contained in the set of preimages of the generic points of the irreducible components of $\cU'_s$ by $spe$ (\cf Theorem \ref{Bo-Lu-Be}). The existence of such $\cU'$ can be deduced for example from Lemma 1.2.11 in \lc Let us still denote by $y_i$, $i=1,\dots,m$ the schematic adherence of points $y_i$ in $\cU'$. The idea is to change $\cU'$ suitably so that $\cU'-\{y_1,\dots,y_m\}$ is affine. Let $D$ be the divisor $\sum_iy_i$ in $\cU'$ and let $\wtilde{D}$ be its reduction in the special fiber $\cU'_s$. Then, by blowing down some of the irreducible components of $\cU'_s$ that are not contained in $Y_s$, where by $Y_s$ we denote the image of $Y$ by the specialization map $\widehat{\cU'}_\eta\to \cU'_s$, and having in mind our initial choice for $\cU'$, we may assume that  $\cU'_s-Y_s$ is just a finite set of closed points. By further blowing down the irreducible components of $\cU'_s$ that are contained in $Y_s$ we may assume that the divisor $\wtilde{D}$ intersects every irreducible component of $\cU'_s$ (in fact, $\cU'_s$ will be the compactification of the canonical reduction of $Y$; furthermore the generic points of the remaining irreducible components of $\cU'_s$ are the images of the points in the Shilov boundary of $Y$ by its canonical reduction, \cf \cite[Proposition 2.4.4.]{Ber90}). In particular it is ample, and by similar arguments as in \cite[Section 7]{Bo-Lu85} one proves as well that $D$ is ample on $\cU'$, that is $\cU$ is affine (\cf \cite[Proof of Proposition 1.2.5]{Bal10}).
\qed

We will use Runge's theorem in the case where $Y$ is an affinoid curve with good reduction and in the next lemma we say to which extent we need to approximate.

\begin{lemma} \label{runge first}
 Let $\vphi$ be an analytic function on a quasi-smooth $k$-affinoid curve with good reduction $Y$ that induces a finite morphism $\vphi:Y\to X=\vphi(Y)\subset\P^1_k$. Then $X$ is an affinoid domain in $\P^1_k$ with good reduction. Furthermore, we can approximate $\vphi$ with a rational function $f$ on $Y'$ that has possible poles only at points $y_1,\dots,y_s$ such that for each $\vt\in TY$, $\sigma(\vphi,\vt)=\sigma(f,\vt)$, $d(\vphi,\vt)=d(f,\vt)$ and finally $\deg(\vphi,\eta)=\deg(f,\eta)$.
\end{lemma}
\proof
The image $X$ is compact and quasi-smooth analytic domain in $\P^1_k$, hence it is an affinoid domain as $Y$ is affinoid. The Proposition 1.5 and comments afterwards in \cite{Col03} imply that $X$ has good reduction as well. 

 Let us fix a small enough prolongation of the morphism $\vphi$, $\vphi':U\to V$. For a point $\vt\in TY$ and the corresponding open annulus/connected component in $A_{\vt}$ in $U-Y$, let $T_{\vt}:A_{\vt}\iso A(0;\rho_{\vt},1)$ be a coordinate on it. Similarly, for $\vv=\wtilde{\vphi}(\vt)$ and the image $A_{\vv}:=\vphi(A_{\vt})$ let $S_{\vv}:A_{\vv}\iso A(0;\rho^{d(\vt)}_{\vt},1)$ be a coordinate. We can write $S_{\vv}=\vphi_{\#}(T_{\vt})=T_{\vt}^{d_{\vt}}\cdot 
h(T_{\vt})$ and $\frac{dS_{\vv}}{dT_{\vt}}=\epsilon_{\vt}\cdot 
T_{\vt}^{\sigma(\vphi,\vt)}\cdot g(T_{\vt} )$. Similarly, if $\eta$ is the Shilov boundary of $Y$ and if $\vt \in T_{\eta}Y$, the corresponding residue class $D_{\vt}$ and its image $D_{\vv}:=\vphi(D_{\vt})$ are open discs, and again after the choice of coordinates $T_{\vt}:D_{\vt}\iso D(0,1^-)$ and $S_{\vv}:=D_{\vv}\iso D(0,1^-)$, we may write $S_{\vv}=\vphi_{\#}(T_{\vt})$. We note that there are only finitely many $\vt\in T_{\eta}Y$  for which $\sigma(\vphi,\vt)\neq 0$, as there are only finitely many critical points for $\vphi$ and as such only finitely many residual classes at $\eta$ that contain them (see Lemma \ref{sigma ram}). For each $\vt\in T_\eta Y \cup TY$, it follows from Lemma \ref{rel annuli}$(i)$ that 
$|\epsilon(\vt)|\geq |d(\vphi,\vt)|$, and since $d(\vphi,\vt)\in\{1,\dots,\deg(\vphi)\}$ 
there exists a global constant $c_1\in\R_{>0}$ such that $|\epsilon(\vt)|>c_1$. On the other hand, for $\vt\in T_{\eta}Y^\dg$ and for $\rho$ close to 1, we 
have the equality 
$|\frac{dS_{\vv}}{dT_{\vt}}|_{\vt,\rho}=|\epsilon|\rho^{\sigma(\vphi,\vt)}$, so as 
long as we fix some $\rho_0\in (0,1)$ that is big enough and only allow 
$\rho\in(\rho_0,1)$, we can bound $|\frac{dS_{\vv}}{dT_{\vt}}|_{\vt,\rho}$ from 
below by some constant $c_2\in\R_{>0}$ not depending on $\vt\in  T_{\eta}Y^\dg$ (only finitely many $\sigma(\vphi,\vt)\neq 0$). In final conclusion, there exists a positive constant $C\in \R_{>0}$ such that for all 
$\vt \in T_{\eta}Y^\dg$ and $\rho$ close to 1, 
$|\frac{dS_{\vv}}{dT_{\vt}}|_{\vt,\rho}>C$. We claim that a rational 
function $f$, such that $|f-\vphi|_{\eta}<\min\{\rho_0C,1\}$ approximates $\vphi$ to the level we need.

Indeed, for $\rho$ close enough to 1 and any $\vt\in T_\eta Y^\dg$, we have (in the second inequality we 
use Lemma \ref{derivative norm}) 
$$|\frac{d}{dT_{\vt}}f-\frac{d}{dT_{\vt}}\vphi|_{\vt,\rho}\leq|\frac{d}{dT_{\vt}
}
|_{\rho}|f-\vphi|_{\vt,\rho}\leq\rho^{-1}\rho_0C\leq 
C<|\frac{d}{dT_{\vt}}\vphi|_{\vt,\rho},
$$
hence $|\frac{d}{dT_{\vt}}f|_{\vt,\rho}=|\frac{d}{dT_{\vt}}\vphi|_{\vt,\rho}$, which 
implies $\sigma(f,\vt)=\sigma(\vphi,\vt)$. Also, as $|\vphi-f|_\eta<1$, $\vphi$ and $f$ induce the same morphism in the reduction, so $d(\vphi,\vt)=d(f,\vt)$ and $\deg(\vphi,\eta)=\deg(f,\eta)$, which finishes the proof.
\qed

\begin{definition}\label{def runge first}
 In the situation of Lemma \ref{runge first}, we say that the rational function $f$ is {\em Runge's first order approximation} of the function $\vphi$. 
\end{definition}
\begin{remark}\label{coh remarks continued} 
Let us keep the notation that we introduced in Remark \ref{coh remarks} $(b)$ and let $X'$ be a simple projectivization of $X$. For each connected component $D_i$, $i=1,\dots,m$, in $X'-X$, let $x_i\in D_i(k)$ be a rational point. Then, the ring of rational functions on $X'-\{x_1,\dots,x_m\}$ is dense in the affinoid algebra $\cO(U_j)$, for every $j\in \N$ as follows from Theorem \ref{runge} ($X'$ is also a simple projectivization of every $U_j$). In particular, the condition $2)$ for the sequence of affinoid domains $U_i$, $i\in \N$, from Remark \ref{coh remarks} $(b)$ is satisfied as well.  
 
\end{remark}

\subsection{The main theorems}

We state the main results of this article. We use notation introduced in \ref{sigma nu}.

\begin{theorem}\label{RH}
 Let $\vphi:Y\to X$ be a finite morphism of quasi-smooth, connected $k$-affinoid curves. Then,
 \begin{equation}\label{RH formula}
 \chi(Y)=\deg(\vphi)\cdot\chi(X)-\sum_{P\in Y(k)}(e_P-1)-\sum_{\vt\in TY}\nu(\vphi,\vt).
 \end{equation}
\end{theorem}
More generaly, we have:
\begin{theorem} \label{RH ft} Let $\vphi:Y\to X$ be a finite morphism of ft $k$-analytic curves. Then
$$
\chi(Y)=\deg(\vphi)\cdot\chi(X)-\sum_{p\in Y(k)}(e_P-1)-\sum_{\vt\in TY}\nu(\vphi,\vt)+\sum_{\vt\in T_{in}Y}\nu(\vphi,\vt).
$$
\end{theorem}
We postpone the proofs until subsection \ref{subsection proofs}.

\subsection{RH formula for affinoid curves in $\P^1_k$}

For the next lemma recall the notions introduced in \ref{sigma nu} and \ref{red morph}.
\begin{lemma}\label{rh rat}
 Let $\Phi:\P^1_k\to \P^1_k$ be a finite morphism (\ie a rational function on $\P^1_k$) and let $x$ be a type two point in $\P^1_k$. Then,
 \beq 
 \sum_{\vt\in T_x\P^1_k} \sigma(\Phi,\vt)=2\deg(\Phi,x)-2.
 \eeq
\end{lemma}
\proof
Let $X$ be an affinoid domain of $\P^1_k$ that has good reduction, $x$ is its Shilov point and such that the restriction $\Phi_{|X}:X\to \Phi(X)=Z$ is a finite morphism. Then, (by definition) $\deg(\Phi,x)=\deg(\Phi_{|X})$. Further, let us fix a coordinate $S$ on the target $\P^1_k$ with respect to which $0\in Z(k)$, while $\infty\in \P^1_k- Z$.  We may also assume that $\vphi(x)$ is sent to the Gauss point with respect to $S$ and let also put $S':=\frac{1}{S}$. 

For $\vt\in T_x\P^1_k$, let $E_{\vt}:=\sum_{P\in D_{\vt}(k)}(e_P-1)$ where  $D_{\vt}$ is the open disc attached to $x$ that corresponds to $\vt$ while $e_P$ is, as usual, the multiplicity of $\Phi$ at the rational point $P$. We note that for any $\vt\in T_xX$, $\sigma(\Phi,\vt)=E_{\vt}$ because of Lemma \ref{sigma ram}. On the other hand, by using valuation polygons it follows that for every $\vt\in TX$, $\sigma(\Phi,\vt)=E_{\vt}-2\sum_{i=1}^{p(\vt)}n_{\vt,i}$ where $n_{\vt,i}$, $i=1,\dots,p(\vt)$ are multiplicities of the poles of the function $S'=\frac{1}{\Phi(T)}$ where $T$ is a coordinate which normalizes disc $D_{\vt}$ that is $ T:D_{\vt}\iso D(0,1^-)$ (note that $\frac{dS'}{dT}=-\frac{\Phi'(T)}{\Phi^2(T)}$). But, the number of poles counted with multiplicity of the function $S'$ is precisely the number of zeros counted with multiplicity in the disc $D_{\vt}$ of the function $S=\Phi(T)$. Summing over $\vt\in TX$ we obtain that $\sum_{\vt\in TX}{\sum_{i=1}^{p(\vt)}n_{\vt,i}}$ is the number of zeros of the function $\Phi$ that are outside of $X$ and in particular this number is equal to $\deg(\Phi)-\deg(\Phi_{|X})=\deg(\Phi)-\deg(\Phi,x)$. To conclude,
\begin{align*}
 \sum_{\vt\in T_x\P^1_k} \sigma(\Phi,\vt)&=\sum_{\vt\in T_x\P^1_k}E_{\vt}+\sum_{\vt\in TX}E_{\vt}-2(\deg(\Phi)-\deg(\Phi,x))\\
 &=2\deg(\Phi)-2-2(\deg(\Phi)-\deg(\Phi,x))\\
 &=2\deg(\Phi,x)-2,
\end{align*}
where we used the classical Riemann-Hurwitz formula $\sum_{\vt\in T_x\P^1_k}E_{\vt}=2\deg(\Phi)-2$. 
\qed
\begin{corollary}\label{RH rat aff good}
 Let $\vphi:Y\to X$ be a finite morphism of $k$-affinoid curves with good reduction and assume that both $Y$ and $X$ are affinoid domains of $\P^1_k$. Then, Riemann-Hurwitz formula \eqref{RH formula} holds for $\vphi$.
\end{corollary}
\proof
 Let us write $\P^1_k-Y=\uplus_{i=1}^\mu D_{Y,i}$ (resp. $\P^1_k-X=\uplus_{j=1}^{\lambda}D_{X,j}$), where $D_{X,i}$'s (resp. $D_{Y,j}$'s) are open discs, and for each $i=1,\dots, \mu$ let $y_i\in D_{Y,i}(k)$ be a rational point. Let $f$ be any first order Runge's approximation of $\vphi$ with respect to the points $y_1,\dots,y_\mu$. Then,  from Lemmas \ref{runge first} and \ref{rh rat} it follows that 
 $$
\sum_{\vt\in T_{\eta}Y^\dg}\sigma(f,\vt)=2\deg(f,\eta)-2=2\deg(\vphi)-2,
$$
where $\eta$ is the Shilov point of $Y$. Also, from Lemma \ref{pullback divisor} we know that $\sum_{\vt\in TY}(d(\vphi,\vt)-1)=\deg(\vphi)\lambda-\mu$, while Theorem \ref{dR groups} implies $\chi(Y)=2-\mu$ and $\chi(X)=2-\lambda$. It follows
$$
\sum_{\vt\in T_{\eta}Y^\dg}\sigma(f,\vt)-\sum_{\vt\in TY}(d(\vphi,\vt)-1)=2\deg(\vphi)-\deg(\vphi)\lambda-2+\mu,
$$
which is, having in mind Lemma \ref{sigma ram} and that $\nu(\vphi,\vt)=\sigma(\vphi,\vt)-d(\vphi,\vt)+1$, just another way to write the RH formula for the morphism $\vphi$. 
\qed
\subsection{Simplifications}

In this subsection we describe a series of simplifications that will allow us reduce the study of RH formula for general morphisms to the case of a morphism where we can apply Theorem \ref{runge}.
\begin{lemma}\label{additivity RH}
 Let $\vphi:Y\to X$ be a finite morphism of ft $k$-affinoid curves, let $\cS$ and $\cT$ be strictly $\vphi$-compatible triangulations of $Y$ and $X$, respectively, such that the restriction of $\vphi$ over the open annuli in $\cA_{\cS}(Y)$ is \'etale (\cf Corollary \ref{consq strict compa} (i)) and assume that $\cS$ has $n>1$ points. For each $s\in \cS$ let $\vphi_s:C_{\cS,s}\to C_{\cT,\vphi(s)}$ be the restriction of $\vphi$ to $C_{\cS,s}$. Then, if RH formula holds for any $n$ out of the $n+1$ morphisms $\vphi:Y\to X$, $\vphi_s:C_{\cS,s}\to C_{\cT,\vphi(s)}$, $s\in \cS$, then it also holds for the remaining one. 
\end{lemma}
\proof
 Suppose that the RH formula holds for all the morphisms $\vphi_s:C_{\cS,s}\to C_{\cT,\vphi(s)}$, that is 
 $$
 \chi(C_{\cS,s})=\deg(\vphi,s)\chi(C_{\cT,\vphi(s)})-\sum_{P\in C_{\cS,s}(k)}(e_P-1)-\sum_{\vt\in TC_{\cS,s}}\nu(\vphi_s,\vt), \quad s\in \cS. 
$$
Note that the classically ramified points are only present in the affinoid domains $C_{\cS,s}$, $s\in \cS$ as $\vphi$ is \'etale over open annuli in $\cA_{\cS}(Y)$. Therefore, summing over $s\in \cS$ and using the Corollary \ref{EP add} we obtain 
$$
\chi(Y)=\deg(\vphi)\chi(X)-\sum_{P\in Y(k)}(e_P-1)-\sum_{\vt\in TY}\nu(\vphi,\vt)-\sum_{s\in \cS}\sum_{\vt\in TC_{\cS,s}-TY}\nu(\vphi_s,\vt).
$$
Furthermore, for each $\vt\in TC_{\cS,s}-TY$ there is an open annulus $A\in \cA_{\cS}(X)$ with an endpoint in $s$ which corresponds to $\vt$. If $s_1\in \cS$ is the other endpoint of $A$, then there exists $\vv\in TC_{S,s_1}-TY$ which also corresponds to $A$ and Lemma \ref{rel annuli} implies that $\nu(\vphi_s,\vt)+\nu(\vphi_{s_1},\vv)=0$, so that 
$$
\sum_{s\in \cS}\sum_{\vt\in TC_{\cS,s}-TY}\nu(\vphi_s,\vt)=0.
$$
In a quite similar way one proceeds also if Riemann-Hurwitz formula holds for $\vphi:Y\to X$ and some $n-1$ out of $n$ morphisms $\vphi_s:C_{\cS,s}:C_{\cT,\vphi(s)}$ to deduce that it also holds for the remaining morphism. 
\qed

The previous lemma with Corollary \ref{RH rat aff good} gives us as a straighforward consequence the following result.
\begin{corollary} \label{RH rat aff} The RH formula holds for finite morphisms between affinoid domains of $\P^1_k$.
\end{corollary}
With a bit more effort one can also prove the following result, for which we omit the proof.
\begin{corollary}\label{general additivity RH}
 Let $\vphi:Y\to X$ be a finite morphism of ft curves, and let $X_i$, $i=1,\dots,n$ be affinoid domains in $X$ such that $X-\uplus_i X_i$ is a disjoint union of open annuli, let $Y_j$, $j=1,\dots,m$ be affinoid domains in $Y$ such that $Y-\uplus Y_j$ is a disjoint union of open annuli, and such that $Y-\uplus_j Y_j=\vphi^{-1}(X-\uplus_i X_i)$ (consequently, for each $j$, $\vphi(Y_j)=X_i$ for some $i$ and the restriction $\vphi_{|Y_i}$ is a finite morphism). Then, if $RH$ holds for $m$ out of the $m+1$ morphisms $\vphi$, $\vphi_{|Y_i}$, $i=1,\dots,m$, it also holds for the remaining one.
\end{corollary}
The next lemma tells us that RH formula is compatible with respect to the composition of finite morphisms. 
\begin{lemma}\label{simplification}
 Let $\vphi:Y\to X$ and $\psi:X\to Z$ be finite morphisms of $k$-affinoid curves. Then, if Riemann-Hurwitz formula holds for two out of the three morphisms $\vphi, \psi, \psi\circ\vphi$, it also holds for the third one. 
\end{lemma}
\proof
Suppose that RH formula holds for morphisms $\psi\circ\vphi$ and $\psi$, the other cases being done in a similar fashion. That is, we have formulas
$$
 \chi(Y)=\deg(\psi\circ\vphi)\chi(Z)-\sum_{P\in Y(k)}(e_P(\psi\circ\vphi)-1)-\sum_{\vt\in TY}\nu(\psi\circ\vphi,\vt)\quad\text{ and }\quad
 $$
 $$
 \chi(X)=\deg(\psi)\chi(Z)-\sum_{P\in X(k)}(e_P(\psi)-1)-\sum_{\vv\in TX}\nu(\psi,\vv).
$$

First we deal with the part concerning classical ramification, so let $D_1,\dots, D_t$ be small enough disjoint strict open discs in $Z$ containing all the branching points of the morphism $\psi\circ\vphi$ (hence also all the branching points of the morphism $\psi$), and so that $\psi^{-1}(D_i)=\uplus_{j=1}^{\alpha(i)}D_{i,j}$ and all $D_{i,j}$ are open discs, and so that $\vphi^{-1}(D_{i,j})=\uplus_{l=1}^{\beta(i,j)}D_{i,j,l}$, and all $D_{i,j,l}$ are open discs as well. Let $\eta_{i,j,l}$ (resp. $\xi_{i,j}$) be the endpoint of $D_{i,j,l}$ (resp. $D_{i,j}$) and let $\vt_{i,j,l}$ (resp. $\vv_{i,j}$) be the tangent point in $T_{\eta_{i,j,l}}Y$ (resp. $T_{\xi_{i,j}}X$) corresponding to the disc $D_{i,j,l}$ (resp. $D_{i,j}$), and for all reasonable $i,j,l$. Then, because of Lemma \ref{sigma ram},
\begin{align*}
\sum_{P\in Y(k)}(e_P(\psi\circ\vphi)-1)&=\sum_{{i,j,l}}\sigma(\psi\circ\vphi,\vt_{i,j,l}), \quad
\sum_{P\in Y(k)}(e_P(\vphi)-1)=\sum_{{i,j,l}}\sigma(\vphi,\vt_{i,j,l}),\\
\sum_{P\in X(k)}(e_P(\psi)-1)&=\sum_{{i,j}}\sigma(\psi,\vv_{i,j}).\\
\end{align*}
Moreover, Lemma \ref{rel annuli} implies that 
\begin{align}
\sum_{i,j,l}\sigma(\psi\circ\vphi,\vt_{i,j,l})&=\sum_{i,j}\sum_l\sigma(\psi\circ\vphi,\vt_{i,j,l})=\sum_{i,j}\big(\deg(\vphi)\sigma(\psi,\vv_{i,j})+\sum_l\sigma(\vphi,\vt_{i,j,l})\big)\nonumber\\
&=\deg(\vphi)\sum_{i,j}\sigma(\psi,\vv_{i,j})+\sum_{i,j,l}\sigma(\vphi,\vt_{i,j,l})\label{composit ram}.
\end{align}

On the other side, using lemmas \ref{rel annuli} and \ref{pullback divisor}, we can write 
\begin{align*}
\sum_{\vt\in TY}\nu(\psi\circ\vphi,\vt)&=\sum_{\vt\in TY}d(\vphi,\vt)\nu(\psi,\wtilde{\vphi}(\vt))+\sum_{\vt\in TY}\nu(\vphi,\vt)
=\sum_{\vv\in TX}\sum_{\substack{\vt\in TY,\\ \wtilde{\vphi}(\vt)=\vv}}d(\vphi,\vt) \nu(\psi,\vv)+\sum_{\vt\in TY}\nu(\vphi,\vt)\\
&=\deg(\vphi)\sum_{\vv\in TX}\nu(\psi,\vv)+\sum_{\vt\in TY}\nu(\vphi,\vt).
\end{align*}
Substituting this, equality \eqref{composit ram} and the RH formula for the morphism $\psi$ in RH formula for the morphism $\psi\circ\vphi$ (and using that $d(\psi\circ\vphi)=d(\psi)d(\vphi)$) we obtain precisely RH formula for the morphism $\vphi:Y\to X$.
\qed

\subsection{Proofs of theorems \ref{RH} and \ref{RH ft}}\label{subsection proofs}

\proof[Proof of Theorem \ref{RH}] First, we may assume that $Y$ and $X$ are with good reduction. Indeed, we may take $\cS$ and $\cT$ to be any strictly compatible triangulations of $Y$ and $X$, respectively, and with $\cS$ having more than one point. Then, if the theorem is true for each morphism $\vphi_s:=\vphi_{|C_{\cS,s}}: C_{\cS,s}\to C_{\cT,\vphi(s)}$, then it is also true for $\vphi:Y\to X$, because of Lemma \ref{additivity RH}.

Second, we may assume that $X$ is an affinoid domain in $\P^1_k$, because of Lemma \ref{simplification} (we can always find a finite morphism from the affinoid $X$ to an affinoid domain in $\P^1_k$).

Let $Y'$ be a simple projectivization of $Y$ and for each connected component $D\in Y'-Y$ we choose a rational point $y_D\in D(k)$. Let $\Phi$ be any first order Runge's approximation of $\vphi$ with respect to points $y_D$. Because $\Phi$ is rational, it induces a finite covering of the projective line $\Phi:Y'\to \P^1$, and one of the connected components of $\Phi^{-1}(X)$ is $Y$. All the other connected components of $\Phi^{-1}(X)$ are isomorphic to affinoid domains in $\P^1_k$, so in particular, Riemann-Hurwitz formula holds for the restriction of $\Phi$ over them (Corollary \ref{RH rat aff}). Let $D_i$, $i\in I$ be the connected components (open discs) of the complement of $X$ in $\P^1_k$, and let $A_i$ be a small enough strict open annulus in $D_i$ that has an endpoint in $X$, such that $\Phi^{-1}(A_i)$ is a disjoint union of open annuli in $Y'$ (thanks to the existence of strictly $\Phi$-compatible triangulations of $Y'$ and $\P^1_k$, we can always find such annuli). Put $Z_i:=D_i-A_i$. All the connected componens of $\Phi^{-1}(Z_i)$, for $i\in I$, are isomorphic to affinoid domains in $\P^1_k$ so in particular RH formula holds for the morphism $\Phi$ restricted over them. Finally, $X$ and $Z_i$, $i\in I$, and connected components of $\Phi^{-1}(X)$ and $\Phi^{-1}(Z_i)$ , $i\in I$, satisfy the conditions of the Corollary \ref{general additivity RH}, and moreover for each component $W\neq Y$ in $\Phi^{-1}(X)\cup_{i\in I} \Phi^{-1}(Z_i)$, RH formula holds for $\Phi_{|W}$. Hence, the same corollary implies the RH formula for the morphism $\vphi$.
\qed

\proof[Proof of Theorem \ref{RH ft}] Let $\cS$ and $\cT$ be any strictly $\vphi$-compatible triangulations of $Y$ and $X$, respectively. Let $\cA_0(Y)$ (resp. $\cA_0(X)$) be the set of open annuli in $\cA_{\cS}(Y)$ (resp. $\cA_{\cT}(X)$) which have only one endpoint in $Y$ (resp. $X$), and put $Y_0:=Y-\uplus_{A\in \cA_0(Y)}$ (resp. $X_0:=X-\uplus_{A\in \cA_0(X)}$). Then, $\vphi_{|Y_0}:Y_0\to X_0$ is a finite morphism of the same degree as $\vphi$ and the RH formula yields 
$$
\chi(Y_0)=\deg(\vphi)\chi(X_0)-\sum_{P\in Y_0(k)}(e_P-1)-\sum_{\vt\in TY_0}\nu(\vphi_{|X_0},\vt).
$$
We note that $\nu(\vphi_{|X_0},\vt)=\nu(\vphi,\vt)$ and that the ramified points in $Y(k)$ are all contained in $Y_0(k)$. Finally, $TY\subseteq TY_0$, and for $\vt\in TY_0-TY$ and the corresponding vector $\vv\in T_{in}Y$,  $\nu(\vphi,\vt)=-\nu(\vphi,\vv)$ because of Lemma \ref{rel annuli}. The proof follows.
\qed

\subsection{Geometric interpretation of the $\sigma$ terms and prolongation of morphisms}

\begin{theorem} \label{sigma general}
Let $\vphi:Y \to X$ be a finite morphism of quasi-smooth $k$-affinoid curves and let $\vphi':Y'\to X'$ be a projectivization of $\vphi$ such that $X'-X$ is a finite union of open discs. Let $\vt \in TY$ and let $U_{\vt}$ be the connected component in $Y'-Y$ which corresponds to $\vt$. Then, $\sigma(\vphi',\vt)=\chi(U_{\vt})-1+\sum_{P\in U_{\vt}(k)}(e_P-1)$.
\end{theorem}
\proof
 The restriction of $\vphi'$ to $U_{\vt}$ is a finite morphism and the image $\vphi'(U_{\vt})$ is one of the connected components of $X'-X$. In particular, $\deg(\vphi'_{|U_{\vt}})=d(\vphi,\vt)$ and RH formula for $\vphi'_{|U_{\vt}}$ yields 
 $$
 \chi(U_{\vt})=d(\vphi',\vt)-\sum_{P\in U_{\vt}(k)}(e_P-1)+\nu(\vphi',\vt)=\sigma(\vphi',\vt)+1-\sum_{P\in U_{\vt}(k)}(e_P-1).\quad\qed
 $$

\begin{remark}
The previous theorem is a generalization of Lemma \ref{sigma ram} and in particular it shows a close relation between $\sigma$ terms and EP characteristic of the corresponding component in the projectivization of the morphism. For example, if $\sigma(\vphi',\vt)<0$ then $U_{\vt}$ cannot be isomorphic to an open disc. In other words, suppose we have a morphism $\vphi:A_1=A(0;r,1)\to A_2=A(0;r^d,1)$ of degree $d$, and let $\sigma(\vphi)$ be the order of the derivative of $\vphi$. Then, if $\sigma(\vphi)<0$ the morphism $\vphi$ does not prolong to a finite morphism of open unit discs. 

It would be interesting to know whether there exist necessary and sufficient conditions one can impose on $\sigma(\vphi)$ so that the morphism $\vphi$ prolongs to a finite morphism of open unit discs.
\end{remark}

\section{Riemann-Hurwitz formula for curves in characteristic {$p>0$}}
 
 We start with a finite morphism $\wtilde{\vphi}': \wtilde{Y'}\to \wtilde{X'}$ of smooth, projective $\tilde{k}$-algebraic curves. When one looks for an RH type of formula for the morphism $\wtilde{\vphi}'$, classically one distingusihes two cases: $\wtilde{\vphi}'$ is separable or $\wtilde{\vphi}'$ is purely inseparable morphism (providing different RH formulas, \cf \cite{Har77,Liu02}). Using the RH formula \eqref{RH formula} for $k$-affinoid curves we provide a method to assign in a uniform way a {\em ramfication divisor} to $\wtilde{\vphi}'$ regardless whether it is separable or purely inseparable.

To begin our argument, let $\wtilde{x}\in 
\wtilde{X'}(\wtilde{k})$ be a rational point, let 
$\wtilde{\vphi}'^{-1}(\wtilde{x})=\{\wtilde{y_1},\dots,\wtilde{y_s}\}$, and let us put $\wtilde{X}:=\wtilde{X'}-  \{\wtilde{x}\}$ and 
$\wtilde{Y}:=\wtilde{Y'}- \{\wtilde{y_1},\dots,\wtilde{y_s}\}$. Then, 
$\wtilde{\vphi}'$ restricts to a finite morphism $\wtilde{\vphi}:\wtilde{Y}\to\wtilde{X}$ 
of smooth {\em affine }curves. 

\begin{lemma}\label{lemma RH p}
 There exists a finite morphism $\vphi:Y\to X$ of quasi-smooth, connected $k$-affinoid curves with good reduction such that:
 \begin{itemize}
  \item[(i)] The reduction of $\vphi$ is $\wtilde{\vphi}$;
  \item[(ii)] $g(Y)=g(\wtilde{Y'})=g(Y')$ and  $g(X)=g(\wtilde{X'})=g(X')$, where $g(*)$ is the genus of the curve "$*$", the genus of an affinoid curve is defined as in \cite[Definition 1.4]{LiuOuverts} and $Y'$ and $X'$ are simple projectivizations of $Y$ and $X$, respectively. 
 \end{itemize}
\end{lemma}
\proof
 To prove $(i)$, let $Y$ and $X$ be any $k$-affinoid curves with reductions $\wtilde{Y}$ and $\wtilde{X}$, respectively, then Theorem \ref{Bo-Lu-Be} implies that $Y$ and $X$ are quasi-smooth and with good reduction. Further, since $\wtilde{X}$ is smooth, \cite[Theorem 1.1]{Col85} implies that there exists a finite morphism $\vphi:Y\to X$ such that the reduction of $\vphi$ is $\wtilde{\vphi}$ (although the result in \lc is presented for the field $\C_p$, it is valid for our field $k$ by the same arguments).
 For $(ii)$, we note that $g(Y)=g(\wtilde{Y}')$ by definition (see \cite[Definition 1.4]{LiuOuverts}), while $g(Y')=g(\wtilde{Y'})$ follows from \cite[Theorem on p. 139]{G-vP80}.
\qed

Any morphism $\vphi:Y\to X$ satisfying the condition $(i)$ in the previous lemma is called a lift of $\wtilde{\vphi}$. We recall that there is a correspondence between the rational points of $\wtilde{Y}'$ and tangent points in $T_{\eta}Y^\dg$, and for a rational point $\tilde{y}\in \wtilde{Y}'(\tilde{k})$, let the corresponding tangent point be denoted by $\vt_{\tilde{y}}$. 

We keep the notation as above.
\begin{theorem}
 The following RH formula for the morphism $\wtilde{\vphi}'$ holds:
 \begin{equation}\label{RH p}
(2-2g(\wtilde{Y'}))=\deg(\wtilde{\vphi})\cdot(2-2g(\wtilde{X'}))-\sum_{\wtilde{y}\in{\wtilde{Y'}(\wtilde{k})}}\sigma(\vphi, \vt_{\wtilde{y}}).
\end{equation}
\end{theorem}
\proof
 The RH formula \eqref{RH formula} for the morphism $\vphi$ yields
$$
\chi(Y)=\deg(\vphi)\chi(X)-\sum_{y\in 
Y(k)}(e_y-1)-\sum_{i=1}^s(\sigma(\vphi,\vt_{\wtilde{y_i}})-d({\vphi,\vt_{\wtilde{y_i}}})
+1),
$$
which is by using Lemma \ref{sigma ram} and Corollary \ref{EP add}, equivalent to 
$$
\chi(Y')=\deg(\vphi)\chi(X')-\sum_{\wtilde{y}\in{\wtilde{Y'}(\wtilde{k})}}
\sigma(\vphi,\vt_{\wtilde{y}}).
$$
The last equality, by Lemma \ref{lemma RH p}, is equivalent to \eqref{RH p}.
\qed

Note that for almost all $\wtilde{y}\in\wtilde{Y'}(\wtilde{k})$, 
$\sigma(\vt_{\wtilde{y}})=0$, so one can define a divisor (which in general depends on the lift $\vphi$)
$$
D(\vphi):=\sum_{\wtilde{y}\in 
\wtilde{Y'}(\wtilde{k})}\sigma(\vphi,\vt_{\wtilde{y}})\wtilde{y}
$$
on $\wtilde{Y}'$ which "measures" the ramification of $\wtilde{\vphi}'$.
When $\wtilde{\vphi}'$ is generically \'etale, the divisor $D(\vphi)$ does not depend on the lift $\vphi$ and is equal to the true ramification divisor of $\wtilde{\vphi}'$, as \cite[Lemma 5.32. 9.]{Bal-Ked} shows. However, the situation is completely different for purely inseparable morphism $\wtilde{\vphi}'$. For example, take 
two morphisms $\vphi_1,\vphi_2:D(0,1)\to D(0,1)$, given by  $S=\vphi_1(T)=T^p$ and $S=\vphi_2(T)=(T-1)^p+1$. 
Then, both $\vphi_1$ and $\vphi_2$ are lifts of the Frobenius covering of the 
affine line by affine line over the field $\wtilde{k}$, while we have 
$\sigma(\vphi_1,\vt_{\wtilde{0}})=p-1\neq0=\sigma(\vphi_2,\vt_{\wtilde{0}})$.

Nevertheless, we propose the following conjecture.
\begin{conjecture}
 For two lifts $\vphi_1$ and $\vphi_2$ of $\wtilde{\vphi}: \wtilde{Y}\to 
\wtilde{X}$, the divisors $D(\vphi_1)$ and $D(\vphi_2)$ are linearly 
equivalent. 
\end{conjecture}

\begin{remark}
 With same kind of reasoning one can provide an RH formula for finite morphisms of $\wtilde{k}$-algebraic curves which are not necessarily smooth provided that there exists a lifting of the morphism (in a suitable sense) to a finite morphism of quasi-smooth $k$-analytic curves.
\end{remark}


\bibliographystyle{plain}
\bibliography{rhformula}

\end{document}